\documentclass[a4paper, 12pt]{article}
\usepackage{enumerate, theorem}
\usepackage{amsmath, amsfonts, amssymb}
\usepackage[height=22.5cm, width=15cm]{geometry}
\usepackage[all]{xy}
\usepackage{mathrsfs, yfonts}

\DeclareMathOperator{\C}{\mathbb{C}}

\newcommand{\A}{\tilde{\mathcal{A}}}

\newcommand{\parag}[1]{\paragraph{\sc{#1.}}}

\newtheorem{thm}{Th\'eor\`eme}[subsection]
\newtheorem{defn}[thm]{D\'efinition}
\newtheorem{cor}[thm]{Corollaire}
\newtheorem{prop}[thm]{Proposition}
\newtheorem{lemma}[thm]{Lemme}

\begin{document}

\title{ Asymptotics of a vanishing period : \\ the quotient themes of a given fresco.}

\author{Daniel Barlet\footnote{Barlet Daniel, Institut Elie Cartan UMR 7502  \newline
Nancy-Universit\'e, CNRS, INRIA  et  Institut Universitaire de France, \newline
BP 239 - F - 54506 Vandoeuvre-l\`es-Nancy Cedex.France. r\newline
e-mail : barlet@iecn.u-nancy.fr}.}

\date{13/01/11}

\maketitle

\section*{Abstract}
In this paper we introduce the word "fresco" to denote a \  $[\lambda]-$primitive monogenic geometric (a,b)-module. The study of this "basic object" (generalized Brieskorn module with one generator) which corresponds to the minimal filtered (regular) differential equation satisfied by a relative de Rham cohomology class, began in [B.09] where the first structure theorems are proved. Then in [B.10] we introduced the notion of theme which corresponds in the \ $[\lambda]-$primitive case  to frescos having a unique Jordan-H{\"o}lder sequence. Themes correspond  to asymptotic expansion of a given vanishing period, so to the image of a fresco in the module of  asymptotic expansions. For a fixed relative de Rham cohomology class (for instance given by  a smooth differential form $d-$closed and $df-$closed) each choice of a vanishing cycle in the spectral eigenspace of the monodromy for the eigenvalue \ $exp(-2i\pi.\lambda)$ \ produces a \ $[\lambda]-$primitive theme, which is a quotient of the fresco associated to the given relative de Rham class  itself. So the problem to determine which theme is a quotient of a given fresco is important to deduce possible asymptotic expansions of the various vanishing period integrals associated to a given relative de Rham class when we change the choice of the vanishing cycle.\\
In the appendix we prove a general existence result which naturally  associate a  fresco to any relative de Rham cohomology class of a proper holomorphic  function of a complex manifold onto a disc.

\parag{AMS Classification} 32 S 25, 32 S 40, 32 S 50.

\parag{Key words} Theme, fresco, (a,b)-module, asymptotic expansion, vanishing period.

\tableofcontents

\section{Introduction}

Let \ $f : X \to D$ \ be an holomorphic function on a connected  complex manifold. Assume that \ $\{df = 0 \} \subset \{ f = 0 \}: = X_0 $. We consider \ $X$ \ as a degenerating family of complex manifolds parametrized by \ $D^* : = D \setminus \{0\}$ \ with a singular member  \ $X_0$ \ at the origin of \ $D$. Let \ $\omega$ \ be a smooth \ $(p+1)-$differential form on \ $X$ \ satisfying \ $d\omega = 0 = df\wedge \omega$. Then in many interesting cases (see for instance [JAG.08] , [B.III] for the case of a function with 1-dimensional singular set and the Appendix for the proper case) the relative family of de Rham cohomology classes induced on the fibers \ $X_s, s \in D^*$ \ of \ $f$ \ by  \ $\omega\big/df$ \ is solution of a minimal filtered differential equation defined from the Gauss-Manin connexion of \ $f$. This object, called a {\bf fresco} is a monogenic regular (a,b)-module satisfying an extra condition, called "geometric", which encodes simultaneously the regularity at \ $0$ \ of the Gauss-Manin connexion, the monodromy theorem and B. Malgrange's positivity theorem.\\
We study the structure of such an object  in order to determine the possible quotient themes of a given fresco. Such a theme corresponds to a possible asymptotic expansion of vanishing periods constructed from \ $\omega$ \ by choosing a vanishing cycle \ $\gamma \in H_p(X_{s_0},\C)$ \ and putting
$$ F_{\gamma}(s) : = \int_{\gamma_s} \ \omega\big/df  $$
where \ $\gamma_s$ \ is the (multivalued) horizontal family of cycles defined from \ $\gamma$ \  in the fibers of \ $f$ \ (see [M.74]).\\
We obtain a rather precise description of the themes of these vanishing periods in term of the structure of the fresco associated to \ $\omega$.\\

We give in the Appendix the existence theorem of the fresco associated to a smooth d-closed and \ $df-$closed form in the case of a proper holomorphic function on a complex manifold. This key result was written in the  preprint  [B.08]  and was not yet published. \\
It is also interesting to indicate to the reader that he may find some "algebraic" explicit computations in [B.09]  of the fresco generated some monomial in the Brieskorn module of the isolated singularity \ $x^5 + y^5 + x^2.y^2 $ \ which is one of the simplest example with a not semi-simple monodromy.

\section{Some known facts.}

\subsection{Regular and geometric (a,b)-modules.}

The main purpose of this paper is to give a precise description of the structure of  a \ $[\lambda]-$primitive monogenic and geometric (a,b)-module ; we shall call such an object a  \ $[\lambda]-$primitive {\bf fresco}. It  corresponds to the \ $[\lambda]-$primitive part of the minimal filtered differential equation satisfied by a relative de Rham cohomology class as indicate in the introduction.\\
 Let us first recall the definition of an (a,b)-module\footnote{For more details on these basic facts, see [B.93]}.

\begin{defn} An (a,b)-module \ $E$ \ is a free finite rank \ $\C[[b]]-$module endowed with an \ $\C-$linear map \ $ a : E \to E$ \ which satisfies the following two conditions :
\begin{itemize}
\item The commutation relation \ $a.b - b.a = b^2$.
\item The map \ $a$ \ is continuous for the \ $b-$adic topology of \ $E$.
\end{itemize}
\end{defn}

Remark that these two conditions imply that for any \ $S \in \C[[b]]$ \ we have
$$ a.S(b) = S(b).a + b^2.S'(b) $$
where \ $S'$ \ is defined via the usual derivation on \ $\C[[b]]$. For a given free rank \ $k$ \ $\C[[b]]-$module with basis \ $e_1, \dots, e_k$, to define a structure of (a,b)-module it is enough to prescribe (arbitrarily) the values of \ $a$ \ on \ $e_1, \dots, e_k$.\\

An alternative way to define  (a,b)-modules is to consider the \ $\C-$algebra 
$$ \A : = \{ \sum_{\nu = 0} ^{\infty} \ P_{\nu}(a).b^{\nu} \} $$
where the \ $P_{\nu}$ \ are polynomials in \ $\C[z]$ \ and where the product by \ $a$ \ is left and right continuous for the \ $b-$adic filtration and satisfies the commutation relation \ $a.b - b.a = b^2$.\\
Then a left \ $\A-$module which is free and finite rank on the subalgebra \ $\C[[b]] \subset \A$ \ is an (a,b)-module and conversely.\\
An (a,b)-module \ $E$ \  has a {\bf simple pole} when we have \ $a.E \subset b.E$ \ and it is {\bf regular} when it is contained in a simple pole (a,b)-module. The regularity is equivalent to the finiteness on \ $\C[[b]]$ \ of the saturation \ $E^{\sharp}$ \ of \ $E$ \  by \ $b^{-1}.a$ \  in \ $E \otimes_{\C[[b]]} \C[[b]][b^{-1}]$.

\begin{itemize}
\item Submodules and quotients of regular (a,b)-modules are regular.
\end{itemize}

Another important property of regular (a,b)-module is the existence of Jordan-H{\"o}lder sequences (J-H. sequences for short).\\
 Recall first that any  regular rank 1 (a,b)-module is characterized up to isomorphism, by a complex number \ $\lambda$ \ and the corresponding isomorphy class is represented by the (a,b)-module  \ $E_{\lambda} : = \C[[b]].e_{\lambda}$ \ where \ $a.e_{\lambda} = \lambda.b.e_{\lambda}$, which is isomorphic to the \ $\A-$module \ $\A\big/\A.(a - \lambda.b)$.\\
Recall also that a submodule \ $F$ \ of the (a,b)-module \ $E$ \ is called  {\bf normal}  when \ $F \cap b.E = b.F$. Normality is a necessary and sufficient condition in order that  the quotient \ $E\big/F$ \ is again an (a,b)-module.\\
A {\bf Jordan-H{\"o}lder sequence} for the rank \ $k$ \ regular (a,b)-module \ $E$ \ is a sequence of normal submodules\footnote{For \ $G \subset F \subset E$ \ submodules with  $F$ \ normal in \ $E$, the normality of \ $G$ \ in \ $F$ \ is equivalent to the normality of \ $G$ \ in \ $E$.}\ $\{0\} = F_0 \subset F_1 \subset \dots F_{k-1} \subset F_k = E$ \ such that the quotients \ $F_j\big/F_{j-1}$ \ for \ $j \in [1,k]$ \ are rank 1 (a,b)-modules. So, to each J-H. sequence of \ $E$, we may associate an ordered sequence of complex numbers \ $\lambda_1, \dots, \lambda_k$ \ such \ $F_j\big/F_{j-1} \simeq E_{\lambda_j}$ \ for each \ $j \in [1,k]$.\\
Existence of J-H. sequence for any regular (a,b)-module and also the following lemma are proved in [B.93].

\begin{lemma}\label{J-H.2}
Let \ $E$ \ be a regular (a,b)-module of rank \ $k$. Up to a permutation,  the set \ $\{exp(-2i\pi.\lambda_j), j \in [1,k]\}$ \ is independant of the choice of the  J-H. sequence of \ $E$. 
Moreover, the sum \ $\sum_{j=1}^k \ \lambda_j$ \ is also independant of the choice of the  J-H. sequence of \ $E$ .
\end{lemma}

The {\bf Bernstein polynomial} of a regular (a,b)-module \ $E$ \ of rank \ $k$ \  is defined as the minimal polynomial of \ $-b^{-1}.a$ \ acting on the \ $k-$dimensional \ $\C-$vector space \ $E^{\sharp}\big/b.E^{\sharp}$. Of course, when \ $E$ \ is the \ $b-$completion of the Brieskorn module of a non constant  germ \ $f : (\C^{n+1}, 0) \to (\C,0)$ \ of holomorphic function with an isolated singularity, we find the "usual" (reduced) Bernstein polynomial of \ $f$ \ (see for instance [K.76] or [Bj.93]).\\

We say that a regular (a,b)-module \ $E$ \ is { \bf geometric} when all roots of its Bernstein polynomial are negative rational numbers. This condition which correspond to M. Kashiwara theorem [K.76],   encodes the monodromy theorem  and the positivity theorem of  B. Malgrange (see [M.75] or the appendix of [B.84] ) extending the situation of  (a,b)-modules deduced from the Gauss-Manin connection of an holomorphic function.\\

Recall that the tensor product of two (a,b)-modules \ $E$ \ and \ $F$ ( see [B.I] ) is defined as the \ $\C[[b]]-$module \ $E\otimes_{\C[[b]]}F$ \ with the \ $\C-$linear endomorphism defined by the rule \ $a.(x\otimes y) = (a.x)\otimes y + x \otimes (a.y)$. The tensor product by a fix (a,b)-module preserves short exact sequences of (a,b)-modules and \ $E_{\lambda}\otimes E_{\mu} \simeq E_{\lambda+\mu}$. So the tensor product of two regular (a,b)-modules is again regular.

\begin{defn}\label{dual}
Let \ $E$ \ be a regular (a,b)-module. The dual \ $E^*$ \ of \ $E$ \ is defined as the \ $\C[[b]]-$module \ $Hom_{\C[[b]]}(E, E_0)$ \ with the \ $\C-$linear map given by
$$ (a.\varphi)(x) =  a.\varphi(x) -  \varphi(a.x)  $$
where \ $E_0 : = \A\big/\A.a \simeq \C[[b]].e_0$ \ with \ $a.e_0 = 0$.
\end{defn}

It is an easy exercice to see that \ $a$ \ acts and satisfies the identity \ $a.b - b.a = b^2$ \ on \ $E^*$ \ with the previous definition. We have \ $E_{\lambda}^* \simeq E_{-\lambda}$ \ and the duality transforms a short exact sequence of (a,b)-modules in a short exact sequence. So the dual of a regular (a,b)-module is again regular. But the dual of a geometric (a,b)-module is almost never geometric. To use duality in the geometric case we shall combine it with tensor product with \ $E_N$ \ where \ $N$ \ is a big enough rational number. Then \ $E^*\otimes E_N$ \ is geometric and \ $(E^*\otimes E_N)^*\otimes E_N \simeq E$. We shall refer to this process as "twisted duality". \\

Define now the left  \ $\A-$module of "formal multivalued expansions"
$$ \Xi : = \oplus_{\lambda \in \mathbb{Q}\cap ]0,1]} \ \Xi_{\lambda} \quad {\rm with} \quad \Xi_{\lambda} : =  \oplus_{j \in \mathbb{N}} \quad \C[[b]].s^{\lambda-1}.\frac{(Log\,s)^j}{j!} $$
with the  action of \ $a$ \ given by
$$ a.(s^{\lambda-1}.\frac{(Log\,s)^j}{j!}) = \lambda.\big[b.s^{\lambda-1}.\frac{(Log\,s)^j}{j!} + b.s^{\lambda-1}.\frac{(Log\,s)^{j-1}}{(j-1)!}\big]$$
for \ $j \geq 1$ \ and \ $a.s^{\lambda-1} = \lambda.b(s^{\lambda-1})$, with, of course, the commutation relations \ $a.S(b) = S(b).a + b^2.S'(b)$ \ for \ $S \in \C[[b]]$.\\
For any geometric (a,b)-module of rank \ $k$,  the vector space \ $Hom_{\A}(E, \Xi)$ \ is of dimension \ $k$ \ and this functor transforms short exact sequences of geometric (a,b)-modules in short exact sequences of finite dimensional vector spaces (see [B.05] for a proof). \\
In the case of the Brieskorn module of an isolated singularity germ of an holomorphic function \ $f$ \  at the origin of \ $\C^{n+1}$ \ this vector space may be identified with the n-th homology group (with complex coefficients) of the Milnor's fiber of \ $f$ (see [B.05]). The correspondance is given by associating to a (vanishing) cycle \ $\gamma$ \ the \ $\A-$linear map
$$ [\omega] \mapsto \big[ \int_{\gamma_s} \ \omega\big/df \big] \in \Xi $$
where \ $\omega \in \Omega^{n+1}_0$,  $\gamma_s$ \ is the multivalued horizontal family of \ $n-$cycles defined by \ $\gamma$ \ in the fibers of \ $f$, and \ $[g]$ \ denotes the formal asymptotic expansion at \ $s = 0$ \ of the multivalued holomorphic function \ $g$.\\

\begin{defn}\label{primitive 1}
A regular (a,b)-module is {\bf \ $[\lambda]-$primitive} (resp. {\bf \ $[\Lambda]-$primitive}), where \ $[\lambda]$ \ is an element (resp. a subset)  in \ $\C\big/\mathbb{Z}$, if all  roots of its Bernstein polynomial are in \ $[- \lambda]$ \ (resp. in \ $[-\Lambda]$).
\end{defn}

\bigskip

If we have a short exact sequence of  (a,b)-modules
$$ 0 \to F \to E \to G \to 0 $$
with \ $E$ \ regular (resp. geometric, resp. $[\lambda]-$primitive) then \ $F$ \ and \ $G$ \ are regular (resp. geometric, resp. $[\lambda]-$primitive). \\
Conversely if \ $F$ \ and \ $G$ \ are regular (resp. geometric, resp. $[\lambda]-$primitive) then \ $E$ \ is regular (resp. geometric, resp. $[\lambda]-$primitive).\\

This implies that \ $E$ \ is \ $[\lambda]-$primitive if and only if it admits a J-H. sequence such that all numbers \ $\lambda_1, \dots, \lambda_k$ \ are in \ $[\lambda]$. And then any J-H. sequence of \ $E$ \ has this property.\\

The following proposition is proved in [B.09] section 1.3.

\begin{prop}\label{primitive 2}
Let \ $E$ \ be a regular (a,b)-module and fix a subset  \ $\Lambda$ \ in \ $ \C\big/\mathbb{Z}$. Then there exists a maximal submodule \ $E[\Lambda]$ \ in \ $E$ \ which is \ $[\Lambda]-$primitive. This submodule is normal in \ $E$.\\
If \ $\{[\lambda_1], \dots, [\lambda_d]\}$ \ is the image in \ $ \C\big/\mathbb{Z}$ \ of the set of the opposite of roots of the Bernstein polynomial of \ $E$, given with an arbitrary order,  there exists a unique sequence \ $0 = F_0 \subset F_1 \subset F_2 \subset F_d = E$ \ of normal submodules of \ $E$ \ such that \ $F_j\big/F_{j-1}$ \ is \ $[\lambda_j]-$primitive for each \ $j \in [1,d]$.
\end{prop}

We  call \ $E[\Lambda]$ \ the {\bf \ $[\Lambda]-$primitive part} of \ $E$.\\

Thanks to this result, to understand what are the possible \ $[\lambda]-$primitive themes which are quotient of a given fresco, it will be enough to work with  \ $[\lambda]-$primitive frescos.\\

Remark also that, in the geometric situation, the choice of a vanishing cycle which belongs to the generalized eigenspace of the monodromy for the eigenvalue \ $exp(-2i\pi.\lambda)$ \ produces vanishing periods with \ $[\lambda]-$primitive  themes.\\

\subsection{Frescos and themes.}

Here we recall some results from [B.09] and [B.10].

\begin{defn}\label{fresco 1}
We shall call {\bf fresco} a geometric (a,b)-module which is generated by one element as an \ $\A-$module.
\end{defn}

\begin{defn}\label{themes 1}
We shall call a {\bf theme}  a fresco which is a sub$-\A-$module of \ $\Xi$.
\end{defn}

Recall that a normal submodule and a quotient by a normal submodule of a fresco (resp. of a theme) is a fresco (resp. is a theme).

A regular rank 1 (a,b)-module is a fresco if and only if it is isomorphic to \ $E_{\lambda}$ \ for some \ $\lambda \in \mathbb{Q}^{+*}$. All rank 1 frescos are themes.
The classification of rank 2 regular (a,b)-modules given in [B.93] gives the list of \ $[\lambda]-$primitive rank 2 frescos  which is the following, where \ $\lambda_1 > 1$ \ is a rational number : 

\begin{align*}
& E =  E \simeq \A\big/\A.(a -\lambda_1.b).(a - (\lambda_1-1).b)  \tag{1}\\
& E \simeq \A\big/\A.(a -\lambda_1.b).(1+ \alpha.b^p)^{-1}.(a - (\lambda_1 + p-1).b) \tag{2}
\end{align*}
where \ $p \in \mathbb{N} \setminus \{0\}$ \ and \ $\alpha \in \C$.\\
The themes in this list are these in \ $(1)$ \ and these in \ $(2)$ \ with \ $\alpha \not= 0 $. For a \ $[\lambda]-$primitive theme  in case \ $(2)$ \ the number \ $\alpha \not= 0$ \ will be called the {\bf parameter} of the theme. \\

For frescos we have a more precise result on the numbers associated to a J-H. sequence :

\begin{prop}\label{fresco 2}
Let \ $E$ \ be a \ $[\lambda]-$primitive fresco and \ $\lambda_1, \dots, \lambda_k$ \ be the numbers associated to a J-H. sequence of \ $E$. Then, up to a permutation, the numbers \ $\lambda_j+j, j\in[1,k]$ \ are independant of the choice of the J-H. sequence.
\end{prop}

The following  structure theorem for frescos will be usefull in the sequel.

\begin{thm}\label{frescos 3}(see [B.09] th.3.4.1)
Let \ $E$ \ be a fresco of rank \ $k$. Then there exists an element \ $P \in \A$ \ which may be written as 
$$ P : = (a - \lambda_1.b)S_1^{-1}.(a - \lambda_2.b)\dots S_{k-1}^{-1}.(a - \lambda_k.b) $$
such that \ $E$ \ is isomorphic to \ $\A\big/\A.P$. Here \ $S_1, \dots, S_{k-1}$ \ are  elements in \ $\C[b]$ \ such that \ $S_j(0) = 1$ \ for each \ $j \in [1,k-1]$.\\
The  element  \ $P_E : = (a - \lambda_1.b)\dots (a - \lambda_k.b) $ \ of \ $\A$ \  is homogeneous in (a,b) and gives the Bernstein polynomial \ $B_E$ \ de \ $E$ \ via the formula
$$(-b)^{-k}.P_E = B_E(-b^{-1}.a) .$$
So \ $P_E \in \A$ \ depends only on the isomorphism class of \ $E$.
\end{thm}

Note that in the case of a fresco the Bernstein polynomial of \ $E$ \ is equal to  the {\bf characteristic polynomial} of the action of \ $-b^{-1}.a$ \ on \ $E^{\sharp}$. This allows a nice formula for a short exact sequence of frescos :

\begin{prop}\label{frescos 4} (see [B.09] prop.3.4.4)
Let \ $0 \to F \to E \to G \to 0$ \ be a short exact sequence of frescos. Then we have the equality in \ $\A$ :
$$ P_E = P_F.P_G \quad {\rm which \ is \  equivalent \ to} \quad  B_E(x) = B_F(x-rk(G)).B_G(x) .$$
\end{prop}

\bigskip

The situation for a \ $[\lambda]-$primitive theme is more rigid :

\begin{prop}\label{themes 2}
A fresco \ $E$ \ is a \ $[\lambda]-$primitive theme if and only if it admits a unique normal rank 1 submodule. In this case the J-H. sequence is unique and contains all normal submodules of \ $E$. The corresponding numbers \ $\lambda_1, \dots, \lambda_k$ \ are such that the sequence \ $\lambda_j+j$ \ is increasing (may-be not strictly).
\end{prop}

\section{Commutation in Jordan-H{\"o}lder sequences.}

In this section we shall study the possible different J-H. sequences of a given \ $[\lambda]-$primitive fresco. Thanks to proposition \ref{primitive 2}, it is easy to see that  the \ $[\lambda]-$primitive assumption does not reduce the generality  of this study.

\subsection{The principal Jordan-H{\"o}lder sequence.}

\begin{defn}\label{Principal J-H.}
Let \ $E$ \ be a \ $[\lambda]-$primitive fresco of rank \ $k$ \ and let  
 $$0 = F_0 \subset F_1 \subset \dots \subset F_k = E$$
 be a J-H. sequence of \ $E$. Then for each \ $j \in [1,k]$ \ we have \ $F_j\big/F_{j-1} \simeq E_{\lambda_j}$, where \ $\lambda_1, \dots, \lambda_k$ \ are in \ $\lambda + \mathbb{N}$.
We shall say that such a J-H. sequence is {\bf principal} when the sequence \  $[1,k]\ni j \mapsto \lambda_j + j$ \ is increasing.
\end{defn}

It is proved in [B.09] prop. 3.5.2  that such a principal J-H. sequence exists for any \ $[\lambda]-$primitive  fresco. Moreover, the corresponding sequence \ $\lambda_1, \dots, \lambda_k$ \ is unique.\\
The following proposition shows much more.

\begin{prop}\label{Uniqueness}
Let \ $E$ \ be a \ $[\lambda]-$primitive fresco. Then its principal J-H. sequence is unique.
\end{prop}

 We shall prove the uniqueness by induction on the rank \ $k$ \ of \ $E$. 

We begin by  the case of rank 2.

\begin{lemma}\label{petit}
Let \ $E$ \ be a rank 2 \ $[\lambda]-$primitive fresco and let \ $\lambda_1,\lambda_2$ \ the numbers corresponding to a principal J-H. sequence of \ $E$ \ (so  \ $\lambda_1+1 \leq \lambda_2+2$). Then the normal rank 1 submodule of \ $E$ \ isomorphic to \ $E_{\lambda_1}$ \ is unique.
\end{lemma}

\parag{Proof} The case \ $\lambda_1+ 1 = \lambda_2+2$ \ is obvious because then \ $E$ \ is a \ $[\lambda]-$primitive theme (see [B.10] corollary 2.1.7). So we may assume that \ $\lambda_2 = \lambda_1 + p_1-1$ \ with \ $p_1 \geq 1$ \ and that  \ $E$ \ is the quotient \ $E \simeq \A\big/\A.(a - \lambda_1.b).(a - \lambda_2.b)$ \ (see the classification of rank 2 frescos in 2.2), because the result is clear when \ $E$ \ is a theme.
We shall use the \ $\C[[b]]-$basis \ $e_1, e_2$ \  of \ $E$ \ where \ $a$ \ is defined by the relations
$$ (a - \lambda_2.b).e_2 = e_1 \quad (a - \lambda_1.b).e_1 = 0. $$
This basis comes from the isomorphism above  \ $E \simeq \A\big/\A.(a - \lambda_1.b).(a - \lambda_2.b)$ \  deduced from the classification of rank 2 frescos with \ $e_2 = [1]$ \ and \ $e_1 = (a - \lambda_2.b).e_2$.\\
Let look for \ $x : = U.e_2 + V.e_1$ \ such that \ $(a - \lambda_1.b).x = 0$. Then we obtain
$$ b^2.U'.e_2 + U.(a - \lambda_2.b).e_2 + (\lambda_2-\lambda_1).b.U.e_2 + b^2.V'.e_1 = 0$$
which is equivalent to the two equations :
$$ b^2.U' + (p_1-1).b.U = 0 \quad {\rm and} \quad U + b^2.V' = 0 .$$
The first equation gives \ $U = 0$ \ for \ $p_1 \geq 2$ \ and \ $U \in \C$ \ for \ $p_1 = 1$. As the second equation implies \ $U(0) = 0$, in all cases \ $U = 0$ \ and \ $V \in \C$. So the solutions are in \ $\C.e_1$. $\hfill \blacksquare$\\

Remark that in the previous lemma, if we assume \ $p_1\geq 1$ \ and \ $E$ \ is not a theme, it may exist infinitely many different normal (rank 1) submodules isomorphic to \ $E_{\lambda_2+1}$. But then, \ $\lambda_2+2 > \lambda_1+1$. See remark 2 following \ref{commute 1} \\

\parag{ proof of proposition \ref{Uniqueness}} As the result is obvious for \ $k = 1$, we may assume \ $k \geq 2$ \ and the result proved in rank \ $\leq k-1$. Let \ $F_j, j \in [1,k]$ \ and \ $G_j, j \in [1,k]$ \ two J-H. principal sequences for \ $E$. As the sequences \ $\lambda_j + j$ \ and \ $\mu_j+j$ \ co{\"i}ncide up to the order and are both increasing, they  co{\"i}ncide. Now let \ $j_0$ \ be the first integer in \ $[1,k]$ \ such that \ $F_{j_0} \not= G_{j_0}$. If \ $j_0 \geq 2$ \ applying the induction hypothesis to \ $E\big/F_{j_0-1}$ \ gives \ $F_{j_0}\big/F_{j_0-1} = G_{j_0}\big/F_{j_0-1}$ \ and so \ $F_{j_0} = G_{j_0}$. \\
So we may assume that \ $j_0 = 1$. Let \ $H$ \ be the normalization\footnote{The smallest normal submodule containing \ $F_1+ G_1$ ;  it has the same rank than \ $F_1+G_1$.} of \ $F_1 + G_1$. As \ $F_1$ \ and \ $G_1$ \ are normal rank 1 and distinct, then \ $H$ \ is a rank 2 normal submodule. It is a \ $[\lambda]-$primitive fresco of rank 2 with two normal rank 1 sub-modules which are isomorphic as \ $\lambda_1 = \mu_1$. Moreover the principal J-H. sequence of \ $H$ \ begins by a normal submodule isomorphic to \ $E_{\lambda_1}$. So the previous lemma  implies \ $F_1 = G_1$. So for any \ $j \in [1,k]$ \ we have \ $F_j = G_j$. $\hfill \blacksquare$\\

\begin{defn}\label{Inv. fond.}
Let \ $E$ \ be a \ $[\lambda]-$primitive fresco and consider a J-H. sequence \ $F_j, j \in [1.k]$ \ of \ $E$. Put \ $F_j\big/F_{j-1} \simeq E_{\lambda_j}$ \ for \ $j \in [1,k]$ \ (with \ $F_0 = \{0\}$). We shall call {\bf fundamental invariants} of \ $E$ \ the (unordered) k-tuple\ $\{ \lambda_j + j, j\in [1,k]\}$.
\end{defn}

Of course this definition makes sens because we know that this (unordered) k-tuple is independant of the choice of the J-H. of \ $E$.\\

Note that when \ $E$ \ is a theme, the uniqueness of the J-H. gives a natural order on this k-tuple. So in the case of a theme the fundamental invariants will be an ordered k-tuple. With this convention, this is compatible with the definition of the fundamental invariants of a \ $[\lambda]-$primitive theme given in [B.10], up to a shift.\\

In the opposite direction, when \ $E$ \ is a semi-simple \ $[\lambda]-$primitive fresco we shall see  (in section 4) that any order of this k-tuple may be realized by a J-H. sequence of \ $E$.\\

\subsection{Commuting in \ $\A$.}

Let \ $E$ \ be a \ $[\lambda]-$primitive fresco. Any isomorphism \ $E \simeq \A\big/\A.P$ \ where \ $P \in \A$ \ is given by
$$ P : = (a - \lambda_1.b).S_1^{-1}.(a - \lambda_2.b) \dots S_{k-1}^{-1}.(a - \lambda_k.b) $$
determines  a J-H. sequence for \ $E$ \ associated to the \ $\C[[b]]-$basis \ $e_1, \dots, e_k$ \ such that  the relations \ $(a - \lambda_j.b).e_j = S_{j-1}.e_{j-1}$ \ hold  for \ $j \in [1,k]$ \ with the convention \ $e_0 = 0$, and the fact that \ $e_k$ \ corresponds, via the prescribed isomorphism, to the class of 1  modulo \ $\A.P$.

\begin{lemma}\label{commute 1}
Let \ $p_1$ \ and \ $p_2$ \ be two positive integers and let \ $\lambda_1 \in \lambda + 2 + \mathbb{N}$, where \ $\lambda \in ]0,1] \cap \mathbb{Q}$. Define \ $P \in \A$ \ as
$$ P : = (a - \lambda_1.b)S_1^{-1}(a - \lambda_2.b)S_2^{-1}(a - \lambda_3.b) $$
where \ $\lambda_{j+1} : = \lambda_j + p_j -1$ \ for \ $j =1,2$, and where \ $S_1, S_2$ \ lie in \ $\C[[b]]$ \ and satisfy \ $S_1(0) = S_2(0) = 1$. We  assume that the coefficient of \ $b^{p_1}$ \ in \ $S_1$ \ vanishes and that the coefficient of \ $b^{p_2}$ \ of \ $S_2$ \ is \ $\alpha \not= 0$.\\
Then if \ $U \in \C[[b]]$ \ is any solution of the differential equation \ $b.U' = p_1.(U - S_1)$, we have
$$P = U^{-1}.(a - (\lambda_2+1).b).(S_1.U^{-2})^{-1}.(a - (\lambda_1-1).b).(U.S_2)^{-1}.(a - \lambda_3.b) .$$
Moreover, there exists an unique choice of \ $U$ \ such that the coefficient of \ $b^{p_1+p_2}$ \ in \ $U.S_2$ \ vanishes.
\end{lemma}

\parag{Proof} The fact that a solution \ $U$ \ of this differential equation satisfies
$$ (a - \lambda_1.b)S_1^{-1}(a - \lambda_2.b) = U^{-1}.(a - (\lambda_2+1).b).S_1^{-1}.U^2.(a - (\lambda_1-1).b).U^{-1} $$
is proved in [B.09] lemma 3.5.1. We use here the case  \ $\delta : = \lambda - \mu = \lambda_2+1 - \lambda_1 = p_1$ \ and the fact that the coefficient of \ $b^{\delta}$ \ in  \ $S_1$ \ vanishes.\\
As the solution \ $U$ \ is unique up to \ $\C.b^{p_1}$, to prove the second assertion let \ $U_0$ \ be the solution with no term in \ $b^{p_1}$. Now the coefficient  \ $\beta(\rho)$ \ of \ $b^{p_1+p_2}$ \ in \ $S_2.U$ \ where  \ $U : = U_0 + \rho.b^{p_1}$, is \ $\beta(\rho) = \beta(0) + \rho.\alpha$. As we assumed that \ $\alpha \not= 0$ \ there exists an unique choice of \ $\rho$ \ for which \ $\beta(\rho) = 0$. $\hfill \blacksquare$\\

\parag{Remarks}
\begin{enumerate}
\item In the situation of the previous lemma the rank 3  fresco \ $E : = \A\big/\A.P$ \ is an extension
$$0 \to E_{\lambda_1} \to E \to T_{\lambda_2,p_2}(\alpha) \to 0 $$
where \ $ T_{\lambda_2,p_2}(\alpha)$ \ is the rank 2  theme with fundamental invariants \ $(\lambda_2, p_2)$ \ and parameter \ $\alpha$\ (see the definition \ref{themes 1}), so 
 $$ T_{\lambda_2,p_2}(\alpha) \simeq \A\big/\A.(a - \lambda_2.b).(1+\alpha.b^{p_2})^{-1}.(a - (\lambda_2+p_2-1).b).$$
\item Let \ $\xi$ \ be in \ $\C^*$ \ and choose \ $\rho : =  (\xi-\beta_0)/\alpha$ \ in the previous proof. Then \ $E$ \ is a extension
$$ 0 \to E_{\lambda_2+1} \to E \to T_{\lambda_1-1,p_1+p_2}(\xi)  \to 0 $$
where \ $T_{\lambda_1-1,p_1+p_2}(\xi) $ \ is the rank 2  theme with fundamental invariants \\ $(\lambda_1-1, p_1+ p_2)$ \ and parameter \ $ \xi$. This shows that we may have infinitely many non isomorphic rank 2 themes as quotients of a given rank 3 \ $[\lambda]-$primitive fresco \ $E$. We have also infinitely many different J-H. sequences with the same quotients: \ $(\lambda_2+1, \lambda_1-1, \lambda_3)$.\\
Note that in this situation we have\footnote{It is easy to see that if \ $Z$ \ is a solution in \ $ \C[[b]]$ \ of the differential equation \ $b.Z' - p_1.Z + S_1 = 0$, then \ $$Ker(a-(\lambda_2+1).b) = \{(r,s) \in \C^2 \ / \ r.(Z.e_1 + b.e_2) + s.b^{p_1}.e_1\}.$$}  \ $\dim_{\C} \big[Ker(a-(\lambda_2+1).b)\big] = 2$.
\end{enumerate}

\begin{cor}\label{commute 2}
In the situation of the previous lemma choose \ $\rho = -\beta(0).\alpha$ \ and denote by \ $V$ \ a solution\footnote{Note that as \ $U.S_2$ \ has no term in \ $b^{p_1+p_2}$ \ with our choice of \ $\rho$, such a solution exists in \ $\C[[b]]$. Moreover \ $V(0) = (U.S_2)(0) = 1$ \ because \ $U(0) = S_1(0) = 1 = S_2(0)$.} of the differential equation \ $b.V' = (p_1+p_2).(V - U.S_2)$. Then \ $P$ \ is equal to
 \begin{equation*}
  U^{-1}.(a - (\lambda_2+1).b).S_1^{-1}.U^2.V^{-1}.(a - (\lambda_3+1).b).(US_2V^{-2})^{-1}.(a- (\lambda_1-2).b).V^{-1} \tag{@}.
  \end{equation*}
 Moreover, the coefficient of \ $b^{p_2}$ \ in \ $S_1.U^{-2}.V$ \ is \ $(p_1+p_2).\alpha\big/p_1 $.
 \end{cor}
 
 \parag{Proof} Of course the choice of \ $\rho$ \ allows to apply again the lemma 3.5.1. of [B.09], with now  \ $\delta = \lambda_3 +1 - (\lambda_1-1) = p_1 + p_2$. This gives \ $(@)$.\\
 Using \ $b.U' = p_1.(U - S_1) $ \ we get
   \begin{align*}
  & b.U'.U^{-2} = p_1.(U^{-1} - S_1.U^{-2}) \quad {\rm and \ with} \quad  Z : = U^{-1} \\
  & b.Z' = -p_1.(Z - S_1.U^{-2}) \quad {\rm and \ then } \\
  & b.Z'.V = -p_1.(Z.V - S_1.U^{-2}.V) \tag{@@}
  \end{align*}
  But using also \ $b.V' = (p_1+p_2).(V - U.S_2)$ \ we get
 $$ b.V'.Z = (p_1+p_2).(V.Z  - S_2) .$$
Adding with \ $(@@)$ \ gives
 $$ b.(V.Z)' - p_2.V.Z = p_1.S_1.U^{-2}.V - (p_1 + p_2).S_2 $$
 which leads to the result, because the left handside has no term in \ $b^{p_2}$. $\hfill \blacksquare$
 
 \bigskip
 
 An obvious consequence of this corollary is that there exists in \ $E$ \ a normal sub-theme isomorphic to \ $T_{\lambda_2+1, p_2}((1+p_2/p_1).\alpha)$, so with fundamental invariants \\
  $(\lambda_2+1,p_2)$ \ and with parameter \ $(1 + p_2/p_1).\alpha$.\\
 Recall that \ $T_{\lambda_2,p_2}(\alpha)$ \ was the rank 2 quotient theme which appears in the principal J-H. of  \ $E$.
 
 \subsection{Some examples.}

We shall give here  some examples, showing the complexity of the non commutative structure of the algebra \ $\A$.

\begin{lemma}\label{exemple 1}
Let \ $x,y,z$ \ non zero complex numbers, \ $\lambda_1$ \ a rational number bigger than 3 and \ $p_1, p_2,p_3$ \ three positive distinct  integers. Assume that \ $p_3$ \ is not a multiple of \ $p_2$, and define \ $\lambda_{j+1} : = \lambda_j+p_j - 1$ \ for \ $j = 1,2,3$. Put 
\begin{enumerate}[]
\item $R_1 : = 1 + x.b^{p_1}, \quad R_3 : = 1+ y.b^{p_3}+ z.b^{p_2+p_3}, \quad  U : = 1 - \frac{z}{y}.b^{p_2}$, 
 \item $ S : = R_1.U \quad {\rm and} \quad T : = U.R_3.$
\end{enumerate}
Then \ $T$ \ has no term in \ $b^{p_2+p_3}$  \ and there exists a solution \ $V \in \C[[b]]$ \ of the differential equation \ $b.V' = (p_2+p_3).(V - T)$.\\
Then the element \ $P : = (a - \lambda_1.b).R_1^{-1}.(a - \lambda_2.b).(a - \lambda_3.b).R_3^{-1}.(a - \lambda_4.b)$ \ in \ $\A$ \ is equal to
\begin{equation*}
 (a - \lambda_1.b).S^{-1}.(a - (\lambda_3+1).b).U^2.V^{-1}.(a - (\lambda_4+1).b).T^{-1}.V^2.(a - (\lambda_2-2).b).V^{-1}. \tag{1}
 \end{equation*}
\end{lemma}

\parag{Proof} A simple computation gives
$$ T = 1 -\frac{z}{y}.b^{p_2} + y.b^{p_3} - \frac{z^2}{y}.b^{2p_2+p_3} $$
and 
$$ V = 1 -\frac{z}{y}\frac{p_2+p_3}{p_3}.b^{p_2} + \frac{p_2+p_3}{p_2}.y.b^{p_3} + \rho.b^{p_2+p_3} +  p_2\frac{z^2}{y}.b^{2p_2+p_3}$$
where \ $\rho$ \ is an arbitrary complex number.\\
Using the lemma 3.5.1. of [B.09] and the fact that \ $U$ satisfies \ $b.U' = p_2.(U - 1)$ \ we get
$$ P = (a-\lambda_1.b).S^{-1}.(a - (\lambda_3+1).b).U^2.(a - (\lambda_2-1).b).T^{-1}.(a - \lambda_4.b).$$
As \ $\lambda_4 = (\lambda_2-1) + p_2+p_3-1$ \ and \ $T$ \ has no term in \ $b^{p_2+p_3}$, we obtain, using again the lemma of {\it loc. cit.}
$$(a - (\lambda_2-1).b).T^{-1}.(a - \lambda_4.b) = V^{-1}.(a-(\lambda_4+1).b)T^{-1}.V^2.(a-(\lambda_2-2).b).V^{-1} $$
if \ $V$ \ is a solution of \ $b.V' = (p_2+p_3).(V - T)$;  this implies \ $(1)$. $\hfill \blacksquare$

\begin{lemma}\label{exemple 1 suite}
In the situation of the previous lemma the rank 4  fresco given by \ 
 $E : = \A\big/\A.P$ \ is not a theme, but we have the following exact sequences :
\begin{align*}
& 0 \to  T_1 \to E \to T_2 \to 0 \\
& 0 \to T_3 \to E \to E_{\lambda_2-2} \to 0.
\end{align*}
where \ $T_1$ \ and \ $T_2$ \ are rank 2 themes and \ $T_3$ \ a rank 3 theme.
\end{lemma}

\parag{Proof} The first exact sequence is consequence of the definition of \ $P$, and the rank 2 theme  \ $T_1$ \ has \ $(\lambda_1,p_1)$ \ as fundamental invariants and \ $x$ \ as parameter ; the rank 2 theme \ $T_2$ \ has \ $(\lambda_3,p_3)$ \ as fundamental invariants and \ $y$ \ as parameter.\\
Let \ $e$ \ be a generator of \ $E$ \ whose annihilator is \ $\A.P$. Then the relation \ $(1)$ \ shows that  \ $\varepsilon : = T^{-1}.V^2.(a - (\lambda_2-2).b).V^{-1}.e$ \ in \ $E$ \ is annihilated by 
 $$Q : =  (a - \lambda_1.b).S^{-1}.(a - (\lambda_3+1).b).U^2.V^{-1}.(a - (\lambda_4+1).b).$$
  So \ $\A.\varepsilon$ \ has rank 3 and is normal because \ $E\big/\A.\varepsilon \simeq E_{\lambda_2-2}$. We shall prove that \ $\A.\varepsilon$ \ is a theme. As \ $\lambda_3+1 = \lambda_1 + p_1 + p_2 -1$ \ and \ $\lambda_4+ 1 = (\lambda_3+1) + p_3 -1$, it is enough to check that the coefficient of\ $b^{p_1+p_2}$ \ in \ $S$ \ and the coefficient in \ $b^{p_3}$ \ in \ $U^{-2}.V$ \ do not vanish. As we have
 \begin{enumerate}[i)]
 \item $ S = 1 + x.b^{p_1} - \frac{z}{y}.b^{p_2} - \frac{xz}{y}.b^{p_1+p_2} $ \\
 \item $ V = 1 -\frac{z}{y}\frac{p_2+p_3}{p_3}.b^{p_2} + \frac{p_2+p_3}{p_2}.y.b^{p_3} + \rho.b^{p_2+p_3} +  p_2\frac{z^2}{y}.b^{2p_2+p_3}$ \\
 \item $U^{-2} = \sum_{n=1}^{\infty} \ n.\big(\frac{z}{y}b^{p_2}\big)^{n-1}$
 \end{enumerate}
these coefficients are respectively equal to \ $-\frac{x.z}{y}$ \ and \ $\frac{p_2+p_3}{p_2}.y$ \ using the fact that \ $p_3$ \ is not a multiple of \ $p_2$. $\hfill \blacksquare$\\

\parag{Remark} Choosing for instance \ $U = 1$ \ gives
$$ P = (a - \lambda_1.b).R_1^{-1}.(a - (\lambda_3+1).b).(a - (\lambda_2-1).b).R_3^{-1}.(a - \lambda_4.b)$$
and then if \ $W$ \ is a  solution of the differential equation \ $b.W' = (p_1+p_2).(W - R_1)$ \  we obtain
$$ P = W^{-1}.(a - (\lambda_3+2).b).(R_1.W^{-2})^{-1}.(a - (\lambda_1-1).b).W^{-1}.(a - (\lambda_2-1).b).R_3^{-1}.(a - \lambda_4.b) $$
and we have an exact sequence
$$ 0 \to E_{\lambda_3+2} \to E \to T_4 \to 0$$
where \ $T_4$ \ is a rank 3 theme and where the corresponding J-H. sequence associated to this exact sequence satifies \ $F_2 = S_1(E)$ \ where \ $S_1(E)$ \  is the maximal  semi-simple  normal submodule of \ $E$ \ (see section 4)  .\\

Note that the first exact sequence corresponds to the principal J-H. of \ $E$. \\
The second gives a J-H. sequence such that its  quotients correspond to the order \\
 $\lambda_1+1, \lambda_3+3, \lambda_4+4, \lambda_2 +2$ \ of the increasing sequence  \ $\lambda_j+j, j\in [1,4]$. The last sequence above corresponds to the order \ $\lambda_3+3, \lambda_1+1,\lambda_2+2, \lambda_4+4$. In this example the semi-simple depth \ $d(E)$ \ of \ $E$ \   (see section 4) is equal to \ $3$.
 
 \parag{Exemple} We give here an example a \ $[\lambda]-$primitive fresco of rank 4 with a J-H. sequence having no non commuting index but which is not semi-simple (see  the section 4 below).\\
 
 Let \ $\lambda_1 > 4$ \ be a rational number and \ $p_1, p_2, p_3$ \ be strictly positive integers. Then consider the fresco
 $$ E : = \A\big/ (a -\lambda_1.b).(a - \lambda_2.b).(1 + b^{p_2+p_3})^{-1}.(a - \lambda_3.b).(a - \lambda_4.b) $$
 where we define \ $\lambda_{j+1} = \lambda_j + p_j -1 $ \ for \ $j =1,2,3$. Then it is clear that all indices of the principal J-H. sequence of \ $E$ \ are commuting indices : for \ $i =1$ \ and \ $i = 3$ \ this is obvious, for \ $i = 2$ \ this results from the commuting lemme 3.5.1 of [B.09] and  the fact that \ $p_2+p_3 > p_3 $ \ as we assume \ $p_2 \geq 1$.
 Now we have the equality in \ $\A$ \ :
 \begin{align*}
 &  (a -\lambda_1.b).(a - \lambda_2.b).(1 + b^{p_1+p_2})^{-1}.(a - \lambda_3.b).(a - \lambda_4.b)  = \\
 & \quad \quad (a -\lambda_1.b).(a - \lambda_2.b).(1 + b^{p_1+p_2})^{-1}.(a - (\lambda_4+1).b).(a - (\lambda_3-1).b) .
 \end{align*}
  This shows, because \ $\lambda_4 + 1 = \lambda_2 + p_2 + p_3 -1$, that there exists a subquotient of rank \ $2$ \ of \ $E$ \ which is a theme;  so \ $E$ \ is not semi-simple. In fact, we produce another J-H. sequence with one non commuting index !

\section{Semi-simple frescos.}

\subsection{The semi-simple filtration.}

\begin{defn}\label{semi-simple}
We shall say that a fresco \ $E$ \ is {\bf semi-simple} if any  quotient of \ $E$ \ which is a \ $[\lambda]-$primitive theme for some \ $[\lambda] \in \mathbb{Q}\big/\mathbb{Z}$ \   is of rank \ $\leq 1$.
\end{defn}

\parag{Remarks}
\begin{enumerate}
\item A  \ $[\lambda]-$primitive theme is semi-simple if and only if it has rank \ $\leq 1$. 
\item An equivalent definition of a  semi-simple fresco is to ask that any \ $\A-$linear map
$$ \varphi : E \to \Xi_{\lambda}$$
for some  \ $[\lambda] \in \mathbb{Q}\big/\mathbb{Z}$ \  has rank \ $\leq 1$. This is a necessary condition because \ $\varphi(E)$ \ is a \ $[\lambda]-$primitive theme which is a quotient of \ $E$. The converse comes from the fact that any  \ $[\lambda]-$primitive theme admits an injective \ $\A-$linear map in \ $\Xi_{\lambda}$. 
\item A fresco is semi-simple if and only if for each \ $[\lambda]$ \ its  \ $[\lambda]-$primitive part (see the propositon \ref{primitive 2}) is semi-simple : if \ $E[\lambda]$ \ is the \ $[\lambda]-$primitive part of \ $E$ \ the restriction map  \ $Hom_{\A}(E, \Xi_{\lambda}) \to Hom_{\A}(E[\lambda], \Xi_{\lambda})$ \ is an isomorphism. For instance,  a theme with only rank \ $\leq  1$ \ \ $[\lambda]-$primitive part for each \ $[\lambda] \in \C\big/\mathbb{Z}$ \  is semi-simple.\\
\end{enumerate}

\begin{lemma}\label{sub and quot}
For any \ $F \subset E$ \ a normal submodule of a semi-simple fresco \ $E$, \ $F$ \ and \ $E\big/F$ \ are semi-simple frescos. 
So any sub-quotient\footnote{By a subquotient \ $H$ \  we mean that there exists \ $G \subset F$ \ normal submodules in \ $E$ \ such that \ $H : = F\big/G$. Remark that \ $H$ \ is a quotient of a normal submodule but also a submodule of a quotient of \ $E$, as \ $F\big/G \subset E\big/G$.} and  of a semi-simple fresco is again a semi-simple fresco.
\end{lemma}

\parag{Proof} As any \ $\A-$linear map \ $\psi : F \to \Xi_{\lambda}$ \ extends to a \ $\A-$linear map \ $\varphi : E \to \Xi_{\lambda}$\ ( see section 2 or [B.05])  the semi-simplicity of \ $E$ \ implies the semi-simplicity of \ $F$. The semi-simplicity of \ $E\big/F$ \ is obvious. $\hfill \blacksquare$

\begin{cor}\label{all J-H}
Let \ $E$ \ be a semi-simple fresco  with rank \ $k$ \ and let \ $\lambda_1, \dots, \lambda_k$ \ be the numbers associated to a J-H. sequence of \ $E$. Let \ $\mu_1, \dots, \mu_k$ \ be a twisted permutation\footnote{This means that the sequence \ $\mu_j+j, j \in [1,k]$ \ is a permutation (in the usual sens) of \ $\lambda_j+j , j\in [1,k]$.} of \ $\lambda_1, \dots, \lambda_k$. Then there exists a J-H. sequence for \ $E$ \ with quotients corresponding to \ $\mu_1, \dots, \mu_k$.
\end{cor}

\parag{Proof} As the symetric group \ $\mathfrak{S}_k$ \  is generated by the transpositions \ $t_{j,j+1}$ \ for \ $j \in [1,k-1]$, it is enough to show that, if \ $E$ \ has a J-H. sequence with quotients given by the numbers \ $\lambda_1, \dots, \lambda_k$ \ then there exists a J-H. sequence for \ $E$ \ with quotients \ $\lambda_1, \dots, \lambda_{j-1}, \lambda_{j+1}+1,\lambda_j-1, \lambda_{j+2}, \dots, \lambda_k$ \ for \ $j \in [1,k-1]$. But \ $G : = F_{j+1}\big/F_{j-1}$ \ is a rank 2  sub-quotient of \ $E$ \ with an exact sequence
$$ 0 \to E_{\lambda_j} \to G \to E_{\lambda_{j+1}} \to 0 .$$
As \ $G$ \ is  a rank 2 semi-simple fresco, it admits also an exact sequence
$$ 0 \to G_1 \to G \to G\big/G_1 \to 0$$
with \ $G_1 \simeq E_{\lambda_{j+1}+1}$ \ and \ $G\big/G_1 \simeq E_{\lambda_j-1}$. Let \ $q : F_{j+1} \to G$ \ be the quotient map. Now the J-H. sequence for \ $E$ \ given by
$$ F_1, \dots, F_{j-1}, q^{-1}(G_1), F_{j+1}, \dots ,F_k = E$$
 satisfies our requirement. $\hfill \blacksquare$
 
  \begin{prop}\label{crit. ss}
 Let \ $E$ \ be a  \ $[\lambda]-$primitive fresco. A necessary and sufficient condition in order that \ $E$ \ is semi-simple is that it admits a J-H. sequence with quotient corresponding to \ $\mu_1, \dots, \mu_k$ \ such that the sequence \ $\mu_j+j$ \ is strictly decreasing.
 \end{prop}
 
 \parag{Remarks}
 \begin{enumerate}
 \item As a fresco is semi-simple if and only if for each \ $[\lambda]$ \ its  \ $[\lambda]-$primitive  part is semi-simple, this proposition gives also a criterium to semi-simplicity for any fresco.
  \item This criterium is a very efficient tool to produce easily examples of semi-simple frescos.
  \end{enumerate}
 
 \parag{Proof} Remark first that if we have, for a  \ $[\lambda]-$primitive fresco \ $E$, a J-H. sequence \ $F_j, j \in [1,k]$ \  such that \ $\lambda_j+j = \lambda_{j+1}+j+1$ \ for some \ $j \in [1,k-1]$, then \ $F_{j+1}\big/F_{j-1}$ \ is a sub-quotient of \ $E$ \ which is a  \ $[\lambda]-$primitive theme of rank  2. So \ $E$ \ is not semi-simple, thanks to the previous corollary. So when a \ $[\lambda]-$primitive fresco \ $E$ \ is semi-simple the principal J-H. sequence corresponds to a strictly increasing sequence \ $\lambda_j+j$. Now, thanks again to the previous corollary we may find a J-H. sequence for \ $E$ \  corresponding to the strictly decreasing order for the sequence \ $\lambda_j+j$. \\

 No let us prove the converse. We shall use the following lemma.
 
 \begin{lemma}\label{utile}
 Let \ $F$ \ be a rank \ $k$ \ semi-simple  \ $[\lambda]-$primitive fresco and let \ $\lambda_j+j$ \ the strictly increasing sequence corresponding to its principal J-H. sequence. Let \ $\mu \in [\lambda]$ \ such that \ $0 < \mu+k+1 < \lambda_1+1$. Then any fresco \ $E$ \ in an exact sequence
 $$ 0 \to F \to E \to E_{\mu} \to 0 $$
 is semi-simple (and \ $[\lambda]-$primitive).
 \end{lemma}
 
 \parag{Proof} Assume that we have a rank 2 quotient \ $\varphi : E \to T$ \ where \ $T$ \ is a \ $[\lambda]-$primitive  theme. Then \ $Ker\, \varphi\, \cap F$ \ is a normal submodule of \ $F$ \ of rank \ $k-2$ \ or \ $k-3$. If \ $Ker\,\varphi \,\cap F$ \ is of rank \ $k-3$, the rank of \ $F\big/Ker\, \varphi\, \cap F$ \ is 2 and it injects in \ $T$ \ via \ $\varphi$. So \ $F\big/Ker\, \varphi \cap F$ \ is a rank 2 \ $[\lambda]-$primitive theme. As it is semi-simple, because \ $F$ \ is semi-simple, we get a contradiction.\\
 So the rank of \ $F\big/Ker\, \varphi\, \cap F$ \ is \ $1$ \ and we have an exact sequence
 $$ 0 \to F\big/Ker\, \varphi\, \cap F \to T \to E\big/F \to 0 .$$
 Put \ $F\big/Ker\, \varphi\, \cap F \simeq E_{\lambda}$. Because \ $T$ \ is a \ $[\lambda]-$primitive  theme, we have the inequality \ $\lambda+1 \leq \mu+2$. But we know that \ $ \lambda_1+1 \leq \lambda+k  $ \  because \ $\lambda + k $ \ is in the set \ $\{ \lambda_j+j, j \in [1,k]\}$ \ and \ $\lambda_1+1$ \ is the infimum of this set. So \ $\lambda_1+1 \leq \mu+k+1 $ \ contradicting our assumption that \ $\mu+k+1 < \lambda_1+1$. $\hfill \blacksquare$
 
 \parag{End of proof of the proposition \ref{crit. ss}} Now we shall prove by induction on the rank of a \ $[\lambda]-$primitive fresco \ $E$ \ that if it admits a J-H. sequence corresponding to a strictly decreasing sequence \ $\mu_j +j$, it is semi-simple. As the result is obvious in rank 1, we may assume \ $k \geq 1$ \ and the result proved for \ $k$. So let \ $E$ \ be a fresco of rank \ $k+1$ \ and let \ $F_j, j\in [1,k+1]$ \ a J-H. sequence for \ $E$ \ corresponding to the strictly decreasing sequence \ $\mu_j+j, j \in [1,k+1]$. Put \ $F_{j}\big/F_{j-1} \simeq E_{\mu_j}$ \ for all \ $j \in [1,k+1]$, define \ $F : = F_k$ \ and \ $\mu : = \mu_{k+1}$; then the induction hypothesis gives that \ $F$ \ is semi-simple and we may apply the previous lemma. $\hfill \blacksquare$\\
 
  \begin{prop}\label{ss base 1}
 Let \ $E$ \ be a  fresco. There exists a unique maximal normal semi-simple submodule \ $S_1(E)$ \ in \ $E$. It contains any (normal) submodule of rank 1  contained  in \ $E$. Moreover, if \ $S_1(E)$ \ is of rank 1, then \ $E$ \ is a \ $[\lambda]-$primitive theme.
 \end{prop}
 
 \parag{Proof} For any \ $\lambda$ \ and any  non zero \ $\varphi \in Hom_{\A}(E, \Xi_{\lambda})$ \ let \ $F_1(\varphi)$ \ be the rank 1 submodule of the \ $[\lambda]-$primitive theme  \ $\varphi(E)$. Now put
 $$ S_1(E) : = \cap_{\lambda}\cap_{\varphi \in Hom_{\A}(E, \Xi_{\lambda})\setminus \{0\}}\quad \big[ \varphi^{-1}(F_1(\varphi))\big] .$$
 Let us prove that \ $S_1(E)$ \ is a normal semi-simple submodule. Normality is obvious as it is an intersection of normal submodules. To prove semi-simplicity, let \ $\psi : S_1(E) \to \Xi_{\lambda}$ \ be a \ $\A-$linear map. Using the surjectivity of the restriction \ $\varphi \in Hom_{\A}(E, \Xi_{\lambda}) \to \varphi_{\vert S_1(E)} \in Hom_{\A}(S_1(E), \Xi_{\lambda})$, (see section 2.1  or  [B.05]), we see immediately that \ $\psi$ \ has  rank \ $\leq 1$. So \  $S_1(E)$ \ is semi-simple. \\
  Now consider a semi-simple normal submodule \ $S$ \ in \ $E$. For any \ $\varphi \in Hom_{\A}(E, \Xi_{\lambda})$ \ the restriction of \ $\varphi$ \ to \ $S$ \ has rank \ $\leq 1$. So \ $\varphi(S)$ \ is contained in the normal rank 1 submodule \ $F_1(\varphi)$ \ of the \ $[\lambda]-$primitive theme \ $\varphi(E)$. So \ $S$ \ is contained in \ $\varphi^{-1}(F_1(\varphi))$ \ for each \ $\varphi$. Then \ $S \subset S_1(E)$, and this proves the maximality of \ $S_1(E)$.\\
  Consider now any rank 1 normal submodule \ $F$ \ of \ $E$. As \ $F$ \ is semi-simple and normal in \ $E$, we have \ $F \subset S_1(E)$. If \ $S_1(E)$ \ is rank 1, there exists an unique rank 1 normal submodule in \ $E$. Then \ $E$ \ is a \ $[\lambda]-$primitive theme, thanks to [B.10] theorem 2.1.6. $\hfill \blacksquare$\\
  
  The following interesting corollary is an obvious consequence of the previous proposition.
  
 \begin{cor}\label{ss base 3}
  Let \ $E$ \ be a   fresco and let \ $\lambda_1, \dots, \lambda_k$ \ be the numbers associated to any J-H. sequence of \ $E$. Let \ $\mu_1, \dots, \mu_d$ \ be the numbers associated to any  J-H. sequence of \ $S_1(E)$. Then, for \ $j \in [1,k]$,  there exists a rank 1 normal submodule of \ $E$ \ isomorphic to \ $E_{\lambda_j+j-1}$ \ if and only if there exists \ $i \in [1,d]$ \ such that we have  \ $\lambda_j+j-1 = \mu_i+i-1$.
  \end{cor}

 Of course, this gives the list of all isomorphy classes of  rank 1 normal submodules  contained in \ $E$. So, using shifted duality, we get also the list of all isomorphy classes of  rank 1 quotients of \ $E$.\\

    \begin{defn}\label{ss base 2}
  Let \ $E$ \ be a  fresco. Define inductively the increasing  sequence \ $S_j(E), j \geq 0$ \ of normal submodules of \ $E$   by putting \ $S_0(E) : = \{0\}$ \ and for \ $j \geq 1$ \\
  \ $ S_j(E)\big/S_{j-1}(E) : = S_1(E\big/S_{j-1}(E))$.  We shall call \ $S_j(E), j\geq 0$ \ the {\bf semi-simple filtration} of \ $E$.  We shall call {\bf semi-simple-depth} of \ $E$ \ ({\bf ss-depth} for short) the first integer \ $d = d(E) \geq 0$ \ such that \ $E = S_d(E)$.
  \end{defn}

   \parag{Example} In the example of lemma \ref{commute 1} let \ $F_2$ \ the second step of the principal J-H. sequence of \ $E$. Then \ $F_2 = S_1(E)$ \ is the maximal  semi-simple normal submodule of \ $E$. This is a consequence of the fact that  \  \ $E$ \ is not semi-simple,\ $F_2$ \  admits a J-H. sequence with  quotients \ $E_{\lambda_2+1}, E_{\lambda_1-1}$ \ with \ $\lambda_2+2 > \lambda_1+1$, so we may apply proposition \ref{crit. ss}. \\

  \begin{prop}\label{ss sequence}
  Let \ $E$ \ be a fresco. Then we have the following properties :
  \begin{enumerate}[i)]
  \item Any \ $[\lambda]-$primitive  sub-theme \ $T$ \ in \ $E$ \ of rank \ $j$ \ is contained in \ $S_j(E)$.
  \item Any  \ $[\lambda]-$primitive quotient theme \ $T$ \ of \ $S_j(E)$ \ has  rank \ $\leq j$. 
  \item For any \ $j \in \mathbb{N}$ \ we have 
  $$ S_j(E) = \cap_{\lambda}\cap_{\varphi \in Hom_{\A}(E, \Xi_{\lambda})} \big[ \varphi^{-1}(F_j(\varphi)) \big]$$
  where \ $F_j(\varphi)$ \ is the normal submodule of rank \ $j$ \ of the \ $[\lambda]-$primitive theme \ $\varphi(E)$, with the convention that \ $F_j(\varphi) = \varphi(E)$ \ when the rank of\ $\varphi$ \ is \ $\leq j$.
  \item The ss-depth of \ $E$ \ is equal to \ $d$ \ if and only if \ $d$ \ is the maximal rank of a \ $[\lambda]-$primitive quotient theme of \ $E$.
  \item The ss-depth of \ $E$ \ is equal to \ $d$ \ if and only if \ $d$ \ is the maximal rank of a normal \ $[\lambda]-$primitive sub-theme of \ $E$.
  \end{enumerate} 
   \end{prop}
   
     \parag{Remarks}
   \begin{enumerate}
   \item By definition of the ss-depth \ $d(E)$ \ of \ $E$ \ the semi-simple filtration is strictly increasing for \ $j \in [0,d(E)]$. 
    \item Let \ $E$ \ be a fresco and \ $N \in \mathbb{Z}$ \ such that \ $E \otimes E_N$ \ is geometric (so is again a fresco). Then \ $E$ \ is semi-simple (resp. a theme) if and only if \ $E \otimes E_N$ \ is semi-simple (resp. a theme). Moreover, in this situation  we have \ $d(E) = d(E \otimes E_N)$.
  \item Let \ $F$ \ be a submodule in a fresco  \ $E$, and denote \ $\tilde{F}$ \ its normalization. Then \  $\tilde{F}$ \ is monogenic (being normal in a monogenic) and geometric.  As there exists \ $N \in \mathbb{N}$ \ such that \ $b^N.\tilde{F} \subset F$,  $\tilde{F}$ \ is a theme for \ $F$ \ a theme. The analog result is also true for a semi-simple \ $F$ :  if \ $\varphi : \tilde{F} \to \Xi_{\lambda}$ \ has rank \ $\geq 2$, as \ $F$ \ has finite codimension in \ $\tilde{F}$, the restriction of \ $\varphi$ \ to \ $F$ \ has also rank \ $\geq 2$ \ which contradicts the semi-simplicity of \ $F$.\\ 
  So we have proved the following two assertions :
  \begin{itemize}
   \item If \ $T \subset E$ \ is a theme in a fresco \ $E$, its normalization is also a theme (of same rank than \ $T$).
   \item If \ $S \subset E$ \ is a semi-simple fresco in a fresco \ $E$, its normalization is also a semi-simple fresco.
 \end{itemize}
 \end{enumerate}
   
   \parag{Proof of proposition \ref{ss sequence}} Let us prove i) by induction on \ $j$. As the case \ $j = 1$ \ is obvious, let us assume that \ $j \geq 2$ \ and that the result is proved for \ $j-1$. Let \ $T$ \ a \ $[\lambda]-$primitive theme in \ $E$, and let \ $F_{j-1}(T)$ \ be its normal submodule of rank \ $j-1$ \ (equal to \ $T$ \ if the rank of \ $T$ \ is less than \ $j-1$). Then by the induction hypothesis, we have \ $F_{j-1}(T) \subset S_{j-1}(E)$. Then we have a \ $\A-$linear map\ $T\big/F_{j-1}(T) \rightarrow E\big/S_{j-1}(E)$. If  the rank of \ $T$ \ is at most \ $j$, then \ $T\big/F_{j-1}(T) $ \ has rank at most 1 and its image is in \ $S_1(E\big/S_{j-1}(E))$. So \ $T \subset S_j(E)$.\\
   To prove ii) we also make an induction on \ $j$. The case \ $j = 1$ \ is obvious. So we may assume \ $j \geq 2$ \ and the result proved for \ $j-1$. Let \ $\varphi : S_j(E )\to T$ \ a surjective map on a \ $[\lambda]-$primitive theme \ $T$. By the inductive hypothesis we have \ $\varphi(S_{j-1}(E)) \subset  F_{j-1}(T)$. So we have an induced surjective map
   $$ \tilde{\varphi} : S_j(E)\big/S_{j-1}(E) \to T\big/F_{j-1}(T).$$
  As \ $S_j(E)\big/S_{j-1}(E)$ \ is semi-simple, the image of \ $\tilde{\varphi}$ \ has rank \ $\leq 1$. It shows that \ $T$ \ has rank \ $\leq j$.\\
  To prove iii) consider first a \ $\A-$linear map\ $\varphi : E \to \Xi_{\lambda}$. As \ $\varphi(E)$ \ is a \ $[\lambda]-$primitive theme, \ $\varphi(S_j(E))$ \ is a \ $[\lambda]-$primitive theme quotient of \ $S_j(E)$. So its rank is \ $\leq j$ \ and we have \ $\varphi(S_j(E)) \subset F_j(\varphi)$.\\
  Conversely, for any \ $\A-$linear map\ $\varphi : E \to \Xi_{\lambda}$, the image \ $\varphi(S_j(E))$ \ is a  \ $[\lambda]-$primitive quotient theme of \ $S_j(E)$. So its rank is \ $\leq j$ \ and it is contained in \ $F_j(\varphi)$.\\
  Let us prove iv). If \ $S_d(E) = E$ \ then any \ $[\lambda]-$primitive sub-theme in \ $E$ \ has rank \ $\leq d$ \ thanks to ii). Conversely, assume that for any \ $[\lambda]$ \  any \ $[\lambda]-$primitive sub-theme of \ $E$ \ has rank \ $\leq d-1$ \ and \ $S_{d-1}(E) \not= E$. Then choose a \ $\A-$linear map \ $\varphi : E \to \Xi_{\lambda}$ \ such that \ $\varphi^{-1}(F_{d-1}(\varphi)) \not= E$. Then \ $\varphi(E)$ \ is a \ $[\lambda]-$primitive theme of rank \ $d$ \ which is a quotient of \ $E$, thanks to the following lemma \ref{clef 0}.\\
  To prove v) let us show that if \ $E$ \ is a fresco and \ $N \gg1 $ \ an integer, then \ $E^*\otimes E_N$ \ is again a fresco and that we have the inequality \ $d(E^*\otimes E_N) \geq d(E)$.\\
  The fact that for \ $N$ \ a large enough integer \ $E^*\otimes E_N$ \ is again a fresco is clear.  Now, as \ $E$ \ has a \ $[\lambda]-$primitive quotient theme of rank \ $d$, then \ $E^*\otimes E_N$ \ has a \ $[-\lambda]-$primitive sub-theme of rank \ $d$.\\
   So we obtain the inequality  \ $d(E^*\otimes E_N) \geq d(E)$ \ from i). Now, using again duality and the fact that \ $[\lambda]-$primitive themes are preserved by \ $\otimes E_N$ \ where \ $N $ \ is a natural  integer, we conclude that  \ $d(E) = d(E^*\otimes E_N)$. Then \ $E$ \ admits a \ $[\lambda]-$primitive normal sub-theme of rank \ $d$.\\
  Conversely, if \ $d$ \ is the maximal rank of a (normal) \ $[\lambda]-$primitive sub-theme of \ $E$, then we have \ $d(E^*\otimes E_N) = d$ \ and \ $d(E) = d$. $\hfill \blacksquare$\\
  
  \begin{lemma}\label{clef 0}
Let \ $E$ \ be a rank \ $k$ \ \ $[\lambda]-$primitive theme and denote by \ $F_j$ \ its normal rank \ $j$ \ submodule. Let \ $x \in E \setminus F_{k-1}$. Then the (a,b)-module \ $\A.x \subset E$ \ is a rank \ $k$ \ theme.
\end{lemma}

\parag{Proof}We  may assume \ $E \subset \Xi_{\lambda}^{(k-1)}$ \ and then (see [B.10]) we have the equality  \ $F_{k-1} = E \cap \Xi_{\lambda}^{(k-2)}$. So \ $x$ \ contains a non zero term with \ $(Log\, s)^{k-1}$ \ and then the result is clear. $\hfill \blacksquare$\\

Our next lemma shows that the semi-simple filtration of a normal submodule of a fresco \ $E$ \ is the trace on this submodule of the semi-simple  filtration of \ $E$.

   \begin{lemma}\label{facile}
    Let \ $E$ \ be a fresco and \ $F$ \ any normal submodule of \ $E$. Then for any \ $j \in \mathbb{N}$ \ we have \ $S_j(E) \cap F = S_j(F)$.
    \end{lemma}
    
    \parag{Proof} By induction on \ $j \geq 1$. First \ $S_1(E) \cap F$ \ is semi-simple in \ $F$ \ so contained in \ $S_1(F)$ \ by definition. But conversely, \ $S_1(F)$ \ is semi-simple, so contained in \ $S_1(E)$ \ and also in \ $F$.\\
    Let assume now  that \ $j \geq 2$ \ and that the result is proved for \ $j-1$. Consider now the quotient \  $E \big/S_{j-1}(E)$. As \ $S_{j-1}(E) \cap F = S_{j-1}(F)$, \ $F\big/S_{j-1}(F)$ \ is  a submodule of \ $E \big/S_{j-1}(E)$. Now by the case \ $j = 1$  \ $S_j(F)\big/S_{j-1}(F)$ \ which is, by definition, \ $S_1(F\big/S_{j-1}(F))$ \ is equal to \ $S_1(E\big/S_{j-1}(E)) \cap (F\big/S_{j-1}(F)) $. So we obtain
    $$ S_j(F)\big/S_{j-1}(F) = (S_j(E)\big/S_{j-1}(E)) \cap (F\big/S_{j-1}(F)) .$$
    This implies the equality \ $S_{j}(F) = S_j(E) \cap F$. $\hfill \blacksquare$\\
    
    \begin{lemma}\label{depth subadditive}
    Let \ $0 \to F \to E \to G \to 0$ \ be a short exact sequence of frescos. Then we have the inequalities
    $$ \sup \{ d(F), d(G)\} \leq d(E) \leq d(F) + d(G) .$$
    \end{lemma}
    
    \parag{Proof} The inequality \ $d(F) \leq d(E)$ \ is obvious from the previous lemma. The inequality \ $d(G) \leq d(E)$ \ is then a consequence of the property iv) in proposition \ref{ss sequence}. \\
    Now let \ $\varphi : E \to \Xi_{\lambda}$ \ be an \ $\A-$linear map with rank \ $\delta$. Then the restriction of \ $\varphi$ \ to \ $F$ \ has rank \ $\leq d(F)$. So \ $\varphi(F)$ \ is contained in \ $T_d$, the normal sub-theme of \ $\varphi(E)$ \ of rank \ $d = d(F)$. The map \ $\tilde{\varphi} : E\big/F \to \Xi_{\lambda}$ \ defined by composition of \ $\varphi$ \ with an injection of the theme \ $\varphi(E)\big/T_d$ \ in \ $\Xi_{\lambda}$ \ has rank \ $\delta - d \leq d(E\big/F)$. So  the inequality \ $\delta \leq d(F) + d(E\big/F)$ \ is proved. $\hfill \blacksquare$\\
    
    \subsection{Co-semi-simple filtration.}
    
    \begin{lemma}\label{duality 1}
    Let \ $E$ \ be a fresco. Then there exists a normal submodule \ $\Sigma^1(E)$ \ which is the minimal normal submodule \ $\Sigma$ \ such that \ $E\big/\Sigma$ \ is semi-simple.
    \end{lemma}
    
    \parag{Proof} First recall that if \ $T$ \ is a theme and \ $T\otimes E_{\delta}$ \ is geometric for some \ $\delta \in \mathbb{Q}$, then \ $T\otimes E_{\delta}$ \ is again a theme.\\
    We shall prove that if \ $E$ \ is a fresco and  if \ $N \in \mathbb{Z}$ \ is such that \ $E \otimes E_N$ \ is again a fresco, we have the equality of submodules in \ $E \otimes E_N$ :
    \begin{equation*}
    S_1(E \otimes E_N) = S_1(E) \otimes E_N .\tag{@}
    \end{equation*}
    As \ $S_1(E)\otimes E_N$ \ is a normal semi-simple submodule of \ $E \otimes E_N$ \ the inclusion \ $\supset$ \ in \ $(@)$ \ is clear.\\
    Conversely, \ $S_1(E \otimes E_N)\otimes E_{-N}$ \ is a semi-simple submodule of \ $E \simeq E\otimes E_N \otimes E_{-N}$. So we obtain \ $S_1(E \otimes E_N)\otimes E_{-N} \subset S_1(E)$ \ and we conclude by tensoring by \ $E_N$.\\
    Now we shall prove that \ $S_1(E^*\otimes E_N)^*\otimes E_N$ \ is a fresco and does not depend of \ $N$, large enough.\\
    Let \ $\lambda_j + j, j \in [1,k]$ \ the  sequence corresponding to  the quotient of a  J-H. of \ $E$. Then let \ $q \in \mathbb{N}$ \ such that \ $\lambda_j+j \in ]k,k+q[$ \ for all \ $j \in [1,k]$. The corresponding J-H. for \ $E^*\otimes E_N$ \ has quotients associated to the numbers \ $-(\lambda_j+j) + k + N $ \ which are in \ $]N-q,N[ \subset ]k,+\infty[$ \ for \ $N > k + q$. So \ $E^*\otimes E_N$ \ is a fresco. Then \ $S_1(E^*\otimes E_N)$ \ is also a fresco and has a J-H. sequence corresponding to numbers in a subset of the previous ones. Dualizing again, we obtain that \ $S_1(E^*\otimes E_N)^*\otimes E_N$ \ has a J-H. sequence with corresponding numbers \ $-(\mu_j+j) + k + N$ \ with \ $\mu_j+j \in ]N-q,N[$. So \ $S_1(E^*\otimes E_N)^*\otimes E_N$ \ is a fresco which is a quotient of \ $E$. \\
    We want to show that this quotient is independant of the choice of \ $N$ \ large enough. This is consequence of the fact that
    \begin{align*}
    & S_1(E\otimes E_{N+1})^* = \big(S_1(E\otimes E_N)\otimes E_1\big)^* \\
    & \qquad  =S_1(E \otimes E_N)^*\otimes E_{-1}
    \end{align*}
    and so
    $$ S_1(E \otimes E_{N+1})^* \otimes E_{N+1} = S_1(E\otimes E_N)^*\otimes E_N .$$
    As \ $S_1(E^*\otimes E_N)$ \ is the maximal semi-simple submodule in \ $E^*\otimes E_N$, we conclude that \ $S_1(E^*\otimes E_N)^*\otimes E_N$ \ is the maximal quotient  of \ $E$ \ which is semi-simple. So we have \ $ \Sigma^1(E) = \Big( E^*\otimes E_N\big/S_1(E^*\otimes E_N)\Big)^*\otimes E_N \subset E . \hfill \blacksquare$\\

        \begin{defn}\label{maximal semi-simple quotients 0}
    Let \ $E$ \ be a fresco and define inductively the normal submodules \ $\Sigma ^j(E)$ \ as follows :
    $\Sigma^0(E) : = E$ \ and \ $\Sigma^{j+1}(E) : = \Sigma^1(\Sigma^j(E))$. We call \ $\Sigma^j(E), j \geq 0$ \ the {\bf co-semi-simple filtration} of \ $E$.
     \end{defn}
     
    Note that \ $\Sigma^j\big/\Sigma^{j+1}$ \ is the maximal semi-simple quotient of \ $\Sigma^j$ \ for each \ $j$.
    
    \begin{lemma}\label{maximal semi-simple quotients 1}
    Let \ $E$ \ be a fresco. The normal submodules \ $\Sigma^j(E)$ \ satisfies the following properties :
    \begin{enumerate}[i)]
    \item For any \ $j \in [0,d(E)-1]$ \ we have \ $\Sigma^{j+1}(E) \subset \Sigma^j(E) \cap S_{d-j-1}(E)$ \ where \ $d : = d(E)$ \ is the ss-depth of \ $E$.
    \item For any \ $[\lambda]-$primitive sub-theme \ $T$ \ of rank \ $t$ \ in \ $\Sigma^j(E)$ \ we have the inclusion \ $F_{t-p}(T) \subset \Sigma^{j+p}$, where \ $p \in [0,t]$ \ and  \ $F_{t-p}(T)$ \ is the rank \ $t-p$ \ normal sub-theme of \ $T$.
    \item Put \ $d : = d(E)$. Then we have \ $d(\Sigma^j(E)) = d - j$ \ for each \ $j \in [0,d]$. This implies that \ $\Sigma^d(E) = \{0\}$ \ and that \ $\Sigma^{d-1}(E) \not= \{0\}$ \  is semi-simple.
    \item For any normal submodule \ $F \subset E$ \ we have \ $ \Sigma^j(F)  \subset \Sigma^j(E) \cap \Sigma^{j-1}(F).$
    \end{enumerate}
    \end{lemma}
    
    \parag{Remarks}
    \begin{enumerate}
    \item The inclusion in i) implies \ $\Sigma^j(E) \subset S_{d-j}(E) \quad \forall j \in [0,d(E)]$.
        \item The filtration \ $\Sigma^j, j \in [0,d]$,  is strictly decreasing because of  iii).   
    \end{enumerate}

    \parag{Proof} Let us prove i) by induction on \ $j \in [0,d(E)-1]$. As i) is obvious for \ $j = 0$ \  assume \ $j \geq 1$ \ and i) proved for \ $j-1$. So we know that \ $\Sigma^j(E) \subset S_{d-j}(E)$. The quotient \ $S_{d-j}(E)\big/S_{d-j-1}(E)$ \ is semi-simple, by definition of \ $S_{d-j}(E)$, and so is its submodule \ $\Sigma^j(E)\big/\Sigma^j(E) \cap S_{d-j-1}(E) $. The definition of \ $\Sigma^{j+1}(E)$ \ implies then that we have \ $\Sigma^{j+1}(E) \subset \Sigma^j(E) \cap S_{d-j-1}(E) $. So i) is proved.\\
    To prove ii) it is enough to show it for \ $p = 1$, by an obvious iteration. By definition \ $\Sigma^j(E)\big/\Sigma^{j+1}(E)$ \ is semi-simple. So is the submodule \ $T\big/T \cap \Sigma^{j+1}(E)$. As it is also a \ $[\lambda]-$primitive theme, its rank is \ $\leq 1$ \ showing that \ $T \cap \Sigma^{j+1}(E)$ \ contains the corank 1 normal submodule \ $F_{t-1}(T)$ \ of \ $T$.\\
    To prove iii) let \ $d : = d(E)$ \ and let \ $T$ \ a sub-theme in \ $E$ \ of rank \ $d$. Then, thanks to ii) with \ $j = 0$, $\Sigma^1(E) \cap T$ \ contains a sub-theme of rank \ $\geq d-1$. So \ $d(\Sigma^1(E)) \geq d-1$. Assume that \ $d(\Sigma^1(E)) = d$, then we obtain, thanks to iteration of the previous inequality, that \ $d(\Sigma^{d-1}(E)) \geq 2$. But from i) we know that \ $ \Sigma^{d-1}(E)) \subset S_1(E)$ \ is semi-simple. This is a contradiction. So we obtain \ $d(\Sigma^1(E)) = d-1$ \ and then \ $d(\Sigma^j(E)) = d - j $ \ for each \ $j \in [0,d]$.\\
    To prove iv) we shall make an induction on \ $j \in [0,d(E)]$. As the case \ $j = 0$ \ is obvious\footnote{with the convention \ $\Sigma^{-1}(G) : = G$.}, assume \ $j \geq 1$ \ and the case \ $j-1$ \ proved. As \ $\Sigma^j(E) \cap F \big/ \Sigma^{j+1}(E)\cap F$ \ is a submodule of \ $\Sigma^j(E)\big/\Sigma^{j+1}(E)$ \ which is semi-simple by definition, it is semi-simple and so we have \ $\Sigma^1(\Sigma^j(E) \cap F) \subset \Sigma^{j+1}(E) \cap F$. Now to conclude, as we know that  \ $\Sigma^j(F) \subset \Sigma^j(E) $,  it is enough to remark that for \ $G \subset H$ \ we have \ $\Sigma^1(G) \subset \Sigma^1(H) \cap G$: as \ $H\big/\Sigma^1(H)$ \ is semi-simple, its submodule \ $G\big/G \cap \Sigma^1(H)$ \ is also semi-simple, and so \ $\Sigma^1(G) $ \ is contained in \ $\Sigma^1(H) \cap G$. $\hfill \blacksquare$
    
    \parag{Remark} We shall prove in section 5 that \ $E\big/S_1(E)$ \ and \ $\Sigma^1(E)$ \ are rank \ $d(E)-1$ \ themes and that any normal rank \ $d(E)$ \ theme in \ $E$ contains \ $\Sigma^1(E)$.

    \subsection{Computation of the ss-depth.}
    
    \begin{defn}\label{non comm. index}
    Let \ $E$ \ be a rank \ $k$ \  \ $[\lambda]-$primitive fresco and consider  \\ $[F] : = \{F_j, j \in [1,k]\}$ \ be any J-H. sequence of \ $E$. We shall say that \ $j \in [1, k-1]$ \ is a {\bf non commuting index} for \ $[F]$ \ if the quotient \ $F_{j+1}\big/F_{j-1}$ \ is a theme. If it is not the case we shall say that \ $j$ \ is a {\bf commuting index}. Note that in this case the quotient \ $F_{j+1}\big/F_{j-1}$ \ is semi-simple.\\
    \end{defn}
    
    \begin{lemma}\label{technique 1}
Let  \ $E$ \ be a rank \ $k \geq 2$ \  \ $[\lambda]-$primitive fresco and let \ $F_j, j \in [1,k]$ \ be a  J-H. sequence of \ $E$. Assume that the \ $\A-$linear map \ $\varphi : F_{k-1} \to \Xi_{\lambda}$ \ has rank \ $\delta \geq 1$ \ and that \ $E\big/ F_{k-2}$ \ is a theme. Then any \ $\tilde{\varphi} : E \to \Xi_{\lambda}$ \ extending \ $\varphi$ \ has rank \ $\delta+1$.
\end{lemma}

\parag{Proof} Let \ $e$ \ be a generator of \ $E$ \ such that \ $(a - \lambda_{k-1}.b).S_{k-1}^{-1}.(a - \lambda_k.b).e$ \ is in \ $F_{k-2}$. So, if \ $\lambda_k = \lambda_{k-1}+p_{k-1} -1$ \ we have either \ $p_{k-1} = 0$ \ or \ $p_{k-1} \geq 1$ \ and the coefficient of \ $b^{p_{k-1}}$ \ in  \ $S_{k-1}$ \ does not vanish. Put \ $\varepsilon : = S_{k-1}^{-1}.(a - \lambda_k.b).e$ ; it  is a generator of \ $F_{k-1}$.  Up to a non zero constant, we may assume that 
 $$\varphi(\varepsilon) -  s^{\lambda_{k-1}-1}.(Log\,s)^{\delta-1}/(\delta-1)! \ \in \ \Xi_{\lambda}^{(\delta-2)} .$$
 Now we want to define \ $\tilde{\varphi}(e) = x$ \ where \ $x$ \ is a solution in \ $\Xi_{\lambda}$ \ of the equation
 $$ (a - \lambda_k.b).x = S_{k-1}.\varphi(e) .$$
 Then it is easy to find that, because the coefficient of \ $b^{p_{k-1}}$ \ in \ $S_{k-1}$ \ is not zero, we have
 $$ x  -  s^{\lambda_k-1}.(Log\,s)^{\delta}/\delta ! \  \in \ \Xi_{\lambda}^{(\delta-1)}.$$
 Now the degree in \ $Log\,s$ \ gives our assertion. $\hfill \blacksquare$\\

\begin{cor}\label{technique 2}
In the situation of the lemma \ref{technique 1} we have the equality  \\
 $d(E) = d(F_{k-1}) + 1$.
\end{cor}

\begin{lemma}\label{pratique}
Let \ $E$ \ be a rank \ $k$ \  \ $[\lambda]-$primitive fresco and \ $[F] : = \{F_j, j \in [1,k]\}$, be any J-H. sequence of \ $E$. Let \ $nci(F)$ \ be the number of {\bf non commuting indices} for the J-H. sequence \ $[F]$. Then we have the inequality
$$ d(E) \geq nci(F) + 1 $$
\end{lemma}

\parag{Proof}  We shall prove this by induction on the rank of \ $E$. The cases of rank 1 and 2 are clear. Let assume \ $k \geq 3$ \ and the inequality proved in rank \ $\leq k-1$. Consider a J-H. sequence \ $F_j, j \in [1,k]$, for \ $E$ \ and assume first that \ $E\big/F_{k-2}$ \ is semi-simple. Then we have, denoting  \ $[G]$ \ the J-H. sequence \ $\{F_j, j \in [1,k-1]\}$, for \ $F_{k-1}$ :
$$ nci(F) = nci(G) \leq d(F_{k-1}) - 1 \leq d(E) - 1 $$
using the induction hypothesis and lemma \ref{depth subadditive} ; it concludes this case.\\
Assume now that \ $E\big/F_{k-2}$ \ is a theme. Then using corollary \ref{technique 2} we have \\ 
$d(E) = d(F_{k-1}) + 1$. So we get using again the inductive hypothesis :
$$ nci(F) = nci(G) + 1 \leq d(F_{k-1}) = d(E) - 1 $$
which concludes the proof. $\hfill \blacksquare$\\

\parag{Remarks} 
\begin{enumerate}
\item This inequality may be strict for several J-H. sequences, including the principal one : there are examples of rank 3  fresco with a principal J-H. sequence  \ $[F]$ \ such that \ $nci(F) = 0$ \ which are not semi-simple (see 3.3 ).
\item We shall see using the corollary of the theorem \ref{struct. fresco 1} (see the remark following \ref{struct. fresco 2}) that for any fresco \ $E$ \  there always exists a J-H. sequence \ $[F]$ \ for which we have the  equality \ $d(E) = nci(F) +1 $.
\end{enumerate}

\subsection{Embedding for a  semi-simple fresco.}

The aim of this paragraph is to prove the following embedding theorem for semi-simple \ $[\lambda]-$primitive frescos.

\begin{prop}\label{embedding 1}
Let \ $E$ \ be a rank \ $k$ \ semi-simple \ $[\lambda]-$primitive fresco. Then there exists an \ $\A-$linear injective map \ $\varphi : E \to  \Xi_{\lambda} \otimes \C^l $ \ if and only if \ $l \geq k$.
\end{prop}

\parag{Proof} To show that the existence of \ $\varphi : E \to \Xi_{\lambda}\otimes \C^l$ \ implies \ $l \geq k$, remark that for any linear form \ $\alpha : \C^l \to \C$ \ the composed map \ $(1\otimes \alpha)\circ\varphi$ \ has  rank  at most \ $1$. So the inequality \ $l \geq k$ \ is clear. To prove that there exists an \ $\A-$linear injective map from \ $E$ \ to \ $\Xi_{\lambda}\otimes \C^k$ \ we shall use the following lemma.

\begin{lemma}\label{embedding 2}
Let \ $\lambda_1, \dots, \lambda_k, k \geq 2,$ be numbers in \ $[\lambda] \in \mathbb{Q}\big/\mathbb{Z}$ \ such that \\
 $\lambda_{j+1} = \lambda_j + p_j -1$ \ for each \ $j \in [1,k-1]$ \ with \ $p_j < 0$. Put
$$ Q : = (a - \lambda_2.b).S_2^{-1} \dots S_{k-1}^{-1}.(a - \lambda_k.b) \quad {\rm and} \quad  P : = (a - \lambda_1.b).S_1^{-1}.Q $$
where \ $S_j, j \in [1,k-1]$ \ are invertible elements in \ $\C[[b]]$. Assume also that \\
 $\lambda_1> k -1$. 
Then there exists an unique element  \ $T \in \C[[b]]$ \ which satisfies
$$ Q.T.s^{\lambda_1-k} = S_1.s^{\lambda_1-1} .$$
Moreover \ $T$ \ is invertible in \ $\C[[b]]$.
\end{lemma}

\parag{Proof} We begin by the proof of the case \ $k = 2$. Then we look for \ $T \in \C[[b]]$ \ such that \ $(a - \lambda_2.b).T.s^{\lambda_1-2} = S_1.s^{\lambda_1-1}$. This equation is equivalent to the differential equation
$$ b.T' - p_1.T = (\lambda_1 -1).S_1 $$
which has an unique solution in \ $\C[[b]]$ \ for any \ $S_1 \in \C[[b]]$ \ because \ $p_1 < 0$. Moreover, we have \ $-p_1.T(0) = (\lambda_1 -1).S_1(0)$, so \ $T$ \ is invertible as  \ $S_1$ \ is invertible.\\
Let now prove the lemma by induction on \ $k \geq 2$. We may assume \ $k \geq 3$ \ and the lemma proved for \ $k-1$. Put \ $Q = (a- \lambda_2.b).S_2^{-1}.R$. Our equation is 
\begin{equation*}
 (a - \lambda_2.b).S_2^{-1}.R.T.s^{\lambda_1- k} = S_1.s^{\lambda_1-1}. \tag{@}
 \end{equation*}
Remark that \ $S_2^{-1}.R.T.s^{\lambda_1-k} = V.s^{\lambda_1-2}$ \ for some \ $V \in \C[[b]]$. So,
let \ $U \in \C[[b]]$ \ the unique solution of the equation
$$ (a - \lambda_2.b).U.s^{\lambda_1-2} = S_1.s^{\lambda_1-1}$$
and consider now the equation in \ $T \in \C[[b]]$ :
\begin{equation*}
 R.T.s^{\lambda_1-k} = S_2.U.s^{\lambda_1-2}. \tag{@@}
 \end{equation*}
The inductive hypothesis shows that there exists an unique invertible \ $T \in \C[[b]]$ \ which is solution of \ $(@@)$. Then it  satisfies \ $(@)$. The uniqueness of the solution  \ $T$ \ of \ $(@)$ \  is consequence of the uniqueness of \ $U$ \ and uniqueness in the inductive hypothesis. $\hfill \blacksquare$

\parag{End of the proof of proposition \ref{embedding 1}} We shall prove the assertion by induction on \ $k$ \ the rank of \ $E$. The case \ $k = 1$ \ is obvious, so we may assume \ $k \geq 2$ \ and that we have a \ $\A-$linear injective map \ $\varphi : F \to \Xi_{\lambda}\otimes \C^{k-1}$ \ for a semi-simple \ $[\lambda]-$primitive fresco \ $F$ \ of rank \ $k-1$. Now using the  notations of the previous lemma and proposition \ref{crit. ss} we may assume that \ $E : = \A\big/\A.P$ \ and put \ $F : = \A\big/\A.Q$. Let \ $e$ \ be a generator of \ $E$ \ with annihilator \ $\A.P$. It induces a generator \ $[e]$ \ of \ $F$ \ with annihilator \ $\A.Q$. Note that \ $F = E\big/E_{\lambda_1}$ \ where \ $E_{\lambda_1} \subset E$ \ is generated by \ $Q.e$. So \ $\varphi([e]) \in  \Xi_{\lambda}\otimes \C^{k-1}$ \ is killed by \ $Q$. Using the previous lemma, we find \ $T \in \C[[b]]$ \ such that \ $Q.T.s^{\lambda_1-k} = S_1.s^{\lambda_1-1}$. So \ $T.s^{\lambda_1-k}$ \ is killed by \ $P = (a - \lambda_1.b).S_1^{-1}.Q$. Define now \ $\psi : E \to \Xi_{\lambda}$ \ via \ $\psi(e) : = T.s^{\lambda_1-k}$ \ and also \ $\Phi : E \to \Xi_{\lambda}\otimes \C^k$ \ as the direct sum of \ $\varphi\circ q$ \ and \ $\psi$, where \ $q : E \to F$ \ is the obvious quotient map given by \ $e \mapsto [e]$. Then \ $\Phi$ \ is injective : if \ $\Phi(x) = 0$ \ then \ $q(x) = 0$ \ so \ $x $ \ is in \ $Ker\, q : = \C[[b]].Q(e) \simeq E_{\lambda_1}$. Then 
 $$\psi(x) = \psi(X.Q.e) = X.Q.T.s^{\lambda_1-k} = X.S_1.(S_1^{-1}.Q.T).s^{\lambda_1-k} = X.S_1.s^{\lambda_1-1} = 0$$ 
if \ $x : = X.Q(e)$ \ with \ $X \in \C[[b]]$. As \ $Q.e$ \ generates a free (rank 1) \ $\C[[b]]-$module and \ $S_1$ \ is invertible in \ $\C[[b]]$,  we conclude that \ $X = 0$ \ and also \ $x= 0$. $\hfill \blacksquare$

\section{Quotient themes of a \ $[\lambda]-$primitive fresco.}

\subsection{Structure theorem for \ $[\lambda]-$primitive frescos.}

Now we are ready to describe the precise  structure of a \ $[\lambda]-$primitive fresco. We shall then  deduce the possible  (\ $[\lambda]-$primitive) quotient themes of  a any given \ $[\lambda]-$primitive fresco. 

\begin{thm}\label{struct. fresco 1}
Let \ $E$ \ be a \ $[\lambda]-$primitive fresco and let \ $d : = d(E)$ \ be its ss-depth. Then there exists  a J-H. sequence \ $G_j, j \in [1,k]$ \ for \ $E$ \ with the following properties :
\begin{enumerate}[i)]
\item The quotient \ $E\big/G_{k-d}$ \ is a theme (with rank \ $d$).
\item We have the equality \ $G_{k-d+1} = S_1(E)$.
\end{enumerate}
\end{thm}

\parag{Remark} By twisted duality we obtain also a J-H. sequence \ $G'_j, j \in [1,k]$, such that \ $G'_{d-1} = \Sigma^1(E)$ \ and \ $G'_d$ \  is a theme ; then \ $E\big/G'_{d-1}$ \ is the maximal semi-simple quotient of \ $E$.\\

\parag{Proof} Let \ $\varphi : E \to \Xi_{\lambda}$ \ an \ $\A-$linear map with rank \ $d$. Denote \ $F_j(\varphi)$ \ the J-H. sequence of the \ $[\lambda]-$primitive theme \ $\varphi(E)$. Put \ $H_j : = \varphi^{-1}(F_j(\varphi))$ \ for \ $j \in [0,d]$. Note that \ $H_0 = Ker\,\varphi$.  We shall show that \ $H_1 = S_1(E)$.\\
To show that \ $H_1$ \ is semi-simple, assume that we have a rank 2 \ $\A-$linear map \ $\theta : H_1 \to \Xi_{\lambda}$. Then \ $E\big/Ker\,\theta$ \ has the following  J-H. sequence
$$ 0 \subset \theta^{-1}(F_1(\theta))\big/Ker\,\theta \subset H_1\big/Ker\,\theta \subset \dots \subset H_d\big/Ker\,\theta .$$
As \ $H_1\big/Ker\,\theta$ \ and \ $(H_{j+1}\big/Ker\,\theta)\Big/(H_{j-1}\big/Ker\,\theta)\simeq H_{j+1}\big/H_{j-1}$ \ are \ $[\lambda]-$primitive rank \ $2$ \  themes for \ $j$ \ in \ $ [1,d-1]$, this would imply that \ $E\big/Ker\,\theta$ \ is a rank \ $d+1$ \  \ $[\lambda]-$primitive theme, contradicting the definition of \ $d$. \\
Then \ $H_1$ \ is semi-simple. But as \ $S_1(E)$ \ is contained in \ $H_1 : = \varphi^{-1}(F_1(\varphi))$ \ the equality \ $H_1 = S_1(E)$ \ is proved.\\
Define now \ $G_{k-d+j} : = H_j$ \ for \ $j \in [0,d]$, and complete the J-H.sequence \\
 $G_j, j \in [1,k]$ \  of \ $E$ \ by choosing a J-H. sequence \ $G_j, j\in [1,k-d]$, for \ $H_0$. $\hfill \blacksquare$\\

We have the following easy consequences of this theorem.

\begin{cor}\label{struct. fresco 2}
Let \ $E$ \ be a \ $[\lambda]-$primitive fresco and let \ $d : = d(E)$ \ be its ss-depth. Then \ $E\big/S_1(E)$ \ and \ $\Sigma^1(E)$ \ are \ $[\lambda]-$primitive themes of rank \ $d-1$ \ and we have
$$ rk(S_1(E)) + d(E) = rk(E) + 1. \quad {\rm and} \quad rk(\Sigma^1(E)) = d(E) -1 .$$
 For each \ $j \in [2,d]$ \ the rank of \ $S_j(E)\big/S_{j-1}(E)$ \ is 1. Moreover, any rank \ $d$ \ quotient theme  \ $T$ \ of \ $E$ \ satisfies \ $E\big/S_1(E) \simeq T\big/F_1(T)$. Dualy, any rank \ $d$ \ theme contains \ $\Sigma^1(E)$.
\end{cor}

\parag{Proof} With the notations of the theorem, let \ $T : = E\big/G_{k-d}$. Then \ $G_{k-d+1}\big/G_{k-d}$ \ is \ $F_1(T)$ \ the unique rank 1 normal submodule of the \ $[\lambda]-$primitive theme \ $T$. Then we obtain that \ $E\big/S_1(E) \simeq E\big/G_{k-d+1} \simeq T\big/F_1(T)$ \ proving our first assertion. The computation of the rank of \ $S_1(E)$ \ follows.\\
Consider now the exact sequence
$$ 0 \to S_1(E) \to E \to T(E) \to 0 $$
where \ $T(E)$ \ is the \ $d-1$ \ $[\lambda]-$primitive theme \ $E\big/S_1(E)$. Dualizing and tensoring by \ $E_N$ \ for \ $N $ \  a large enough integer gives the exact sequence
$$ 0 \to T(E)^*\otimes E_N  \to E^*\otimes E_N \to S_1(E)^*\otimes E_N \to 0 $$
where \ $S_1(E)^*\otimes E_N$ \ is semi-simple. This implies that \ $\Sigma^1(E^*\otimes E_N) \subset T(E)^*\otimes E_N$. And we have equality because we know that \ $E^*\otimes E_N$ \ contains a rank  \ $d$ \ theme, so the dimension of \ $\Sigma^1(E^*\otimes E_N)$ \ is at most \ $d(E)-1$. Then we conclude that \ $\Sigma^1(E)$ \ is a theme and that its rank   is \ $d(E)-1$.\\
We know that \ $E\big/S_1(E)$ \ is a theme, so \ $S_1(E\big/S_1(E)) = F_1(E\big/S_1(E))$ \ is rank \ $1$ \ for \ $d(E) \geq 2$. A similar arguement shows that \ $S_j(E)\big/S_{j-1}(E)$ \ for each \ $j \in [2,d]$. In fact it is naturally isomorphic to \ $F_{j-1}(T)\big/F_{j-2}(T)$ \ where \ $T : = T(E)$.\\
Let \ $T'$ \ any rank \ $d$ \ quotient theme of \ $E$. As \ $E$ \ is \ $[\lambda]-$primitive, so is \ $T'$ \ and we may assume that \ $T' : = \varphi(E)$ \ where \ $\varphi \in Hom_{\A}(E,\Xi_{\lambda})$. But now \ $S_1(E)$ \ is in \ $\varphi^{-1}(F_1(T'))$, so we have a surjective map, induced by \ $ \varphi$ :
$$\tilde{\varphi} : E\big/S_1(E) \to T'\big/F_1(T') $$
between two \ $[\lambda]-$primitive themes of the same rank \ $d-1$.
 This must be an isomorphism. $\hfill \blacksquare$
 
 \parag{Remark} Building a J-H. sequence of \ $E$ \ via the exact sequence
 $$ 0 \to S_1(E) \to E \to T(E) \to 0 $$
 we find that the number of non commuting indices in such a J-H. sequence is exactely \ $d(E) - 1$. Of course we may put the non commuting indices at the beginning by using the exact sequence
 $$ 0 \Sigma^1(E) \to E \to  E\big/\Sigma^1(E) \to 0 .$$

 \subsection{Embedding dimension for a \ $[\lambda]-$primitive fresco.}
 
 From the embedding result in the semi-simple case \ref{embedding 1} and the structure theorem \ref{struct. fresco 1} we shall deduce precise embedding theorem for \ $[\lambda]-$primitive frescos.

\begin{prop}\label{embedding 3}
Let \ $E$ \ be a \ $[\lambda]-$primitive fresco. Then there exists an injective \ $\A-$linear map \ $\varphi : E \to \Xi_{\lambda}\otimes \C^l$ \ if and only if \ $l \geq rk(E) - d(E) + 1$.
\end{prop}

\parag{Proof} If we have such a \ $\varphi$ \ its restriction to \ $S_1(E)$ \ is an embedding and so \ $l \geq rk(S_1(E)) = rk(E) - d(E) + 1$.\\
Conversely, we shall prove that there exists and embedding of \ $E$ \ in \ $\Xi_{\lambda} \otimes \C^l$ \ with \ $l = rk(E) - d(E) +1$. By the proposition  \ref{embedding 1} we may begin with an injective \ $\A-$linear map \ $\varphi : S_1(E) \to \Xi_{\lambda}\otimes \C^l$ \ with \ $l : = rk(E) - d(E) + 1$. Put \ $\varphi : = \oplus_{i=1}^l \ \varphi_i $ \ where \ $\varphi_i : S_1(E) \to \Xi_{\lambda}$. Now, for each \ $i \in [1,l]$ \ we may find an extension \ $\Phi_i : E \to \Xi_{\lambda}$ \ to \ $\varphi_i$ \ thanks to the surjectivity of \ $Hom_{\A}(E, \Xi_{\lambda}) \to Hom_{\A}(S_1(E), \Xi_{\lambda})$ \ (see section 2.1).Then we may define  \ $\Phi : = \oplus_{i\in [1,l]} \Phi_i : E \to \Xi_{\lambda}\otimes \C^l$ \ which is an extension of \ $\varphi$ \ to \ $E$. Moreover, such an extension is injective because its kernel cannot meet non trivialy \ $S_1(E)$ \ and so does not contain any normal rank 1 submodule of \ $E$. So the kernel has to be \ $\{0\}$. $\hfill \blacksquare$

\bigskip

\subsection{Quotient themes of a \ $[\lambda]-$primitive fresco.}

Now we shall describe all quotient themes of a given \ $[\lambda]-$primitive fresco. We begin by the description of quotient themes of maximal rank.

\begin{prop}\label{maximal quotient themes}
Let \ $E$ \ be a \ $[\lambda]-$primitive fresco of rank \ $k$. Then any rank \ $d : = d(E)$ \ quotient theme of \ $E$ \ is obtained as follows : let \ $K$ \ be a  corank \ $1$ \  normal submodule of \ $S_1(E)$ \ and assume that \ $K \cap L(E) = \{0\}$ \ where \ $L(E) : = \Sigma^1(S_2(E))$. Then \ $E\big/ K$ \ is a rank \ $d$ \ theme.
\end{prop}

\parag{Proof} Consider \ $K \subset S_1(E)$ \ a corank 1 normal submodule in \ $S_1(E)$ \ such that \ $K \cap L(E) = \{0\}$. By definition of \ $L(E)$ \ the quotient \ $S_2(E)\big/K$ \ is rank 2 and not semi-simple. So it is a theme. We have the following Jordan-H{\"o}lder sequence for \ $E\big/K$ :
\begin{equation*}
0 \subset  S_1(E)\big/K \subset S_2(E)\big/K \subset \dots \subset S_d(E)\big/K = E\big/ K \tag{@} 
\end{equation*}
But \ $S_2(E)\big/K$ \ and each \ $S_{j+2}(E)\big/S_j(E)$ \ for \ $j \in [1,d-2]$ \ is a theme of rank 2. So from [B.10] we conclude that \ $E\big/K$ \ is a rank d  theme.\\
Conversely, if \ $E\big/K$ \ is a rank \ $d$ \ theme, consider \ $S_1(E)\big/S_1(E) \cap K \hookrightarrow E\big/K$. As \ $S_1(E)\big/S_1(E) \cap K$ \ is semi-simple and \ $E\big/K$ \ is a theme, the rank of \ $S_1(E)\big/S_1(E) \cap K$ \ is at most 1. It is not \ $0$ \ because \ $K$ \ has rank \ $k-d$ \ and \ $S_1(E)$ \ has rank \ $k-d+1$.  So \ $K$ \ is contained in \ $S_1(E)$ \ and has corank 1 in it.  If \ $K $ \ contains \ $L(E)$ \ then \ $S_2(E)\big/K$ \ is semi-simple and of rank 2, if we assume \ $d \geq 2$. But it is contained in \ $E\big/K$ \ which is a theme, so we get a contradiction.\\
 For \ $d = 1$ \ in the previous proposition (so \ $E$ \ semi-simple) we have \ $S_2(E) = E$ \ so \ $L(E) = \{0\}$ \ and any corank 1 normal submodule of \ $S_1(E) = E$ \ gives a rank 1 quotient which is, of course a rank 1 theme. $\hfill \blacksquare$\\
 
 In the statement of the theorem we shall denote by \ $L_j$ \ for \ $j \in [1,d-1]$ \ the rank \ $j$ \ theme defined as \ $\Sigma^1(S_{j+1}(E))$. So, by definition, \ $S_{j+1}(E)\big/\Sigma_j$ \ is semi-simple, and \ $L_j$ \ is a normal submodule which is minimal for this property. We have seen that \ $L_j$ \ is then a theme with rank \ $d(S_{j+1}(E)) -1$ ; as we know that \ $d(S_{j+1}(E)) = j+1$, the rank of \ $L_j$ \ is \ $j$. In fact we have \ $L_j = F_j(\Sigma^1(E))$ \ for each \ $j \in [0,d-1]$ \ and so \ $L_1 = L(E)$.

\begin{thm}\label{quotient themes}
Let \ $E$ \ be a \ $[\lambda]-$primitive fresco of rank \ $k$. Put \ $d : = d(E)$ \ and denote by \ $S_j, j \in [1.d]$ \ the semi-simple filtration of \ $E$. Assume that \ $d \geq 2$ \ and let \ $\Sigma_j$ \ be the first term of the co-semi-simple filtration of \ $S_{j+1}$.\\
Let \ $K \subset E$ \ be a normal submodule such that \ $E\big/K$ \ is a theme. Then we have the following possibilities :
\begin{enumerate}
\item If \ $K$ \ contains \ $S_1(E)$, then \ $E\big/K$ \ is a quotient theme of the rank \ $d-1$ \ theme \ $E\big/S_1(E)$ \  and we have exactlely one such quotient for each rank in \ $[1, d-1]$.
\item If \ $K \cap L_1 = \{0\}$, then \ $E\big/K$ \ is a quotient of the rank \ $d$ \ theme   \ $E\big/K \cap S_1(E)$ \ which belongs to the quotient themes described in the previous proposition.
\item   If \ $K $ \ contains \ $L_{j_0}$ \ but not \ $L_{j_0+1}$ \ (so \ $K \cap \Sigma^1(E)$ \ has rank \ $j_0$),  we may apply the previous case to \ $E' : = E\big/L_{j_0}$ \ and \ $K' : = K\big/L_{j_0}$ \ and find that \ $E\big/K = E'\big/K'$ \ is a quotient of the rank \ $d-j_0 $ \ quotient theme \ $E'\big/K' \cap S_1(E')$. In this situation we have  \ $ S_1(E') = S_{j_0+1}(E)\big/L_{j_0}$ \ and \ $d(E') = d- j_0$, with \ $ rk(E') = rk(E) - j_0$. So the rank of \ $E\big/K$ \ is at most \ $d - j_0$.
\end{enumerate}
\end{thm}

\parag{Remarks}
\begin{enumerate}[i)]
\item The case 2 of the previous theorem is the case 3 with \ $j_0 = 0$. We emphasis on this case because the rank d quotients themes is the most interesting case.
\item Let \ $\varphi : E \to T$ \ be a surjective \ $\A-$linear map on a rank \ $\delta \geq 2$ \ theme \ $T$. Then \ $S_1(E) \subset \varphi^{-1}(F_1(T))$ \ and so the map \ $\varphi$ \ induces a surjection \ $E\big/S_1(E) \to T\big/F_1(T)$. As \ $E\big/S_1(E)$ \ is a \ $[\lambda]-$primitive  theme, it has an unique quotient of rank \ $\delta - 1$. So the quotient theme \ $T\big/F_1(T)$ \ depends only on \ $\delta$ \ and \ $E$, not on \ $T$. In the case \ $\delta = d(E)$ \ we find that \ $T\big/F_1(T) \simeq E\big/S_1(E)$ \ for any choice of \ $T$.
\item The previous theorem gives very few information on the rank 1 quotients, because any corank \ $1$ \ normal submodule contains \ $\Sigma^!(E)$. They will be described in the next proposition.
\end{enumerate}

\parag{Proof} The first case is clear.\\
Assume that \ $L_1 \cap K = \{0\}$;  then \ $K$ \ does not contain \ $S_1$. But \ $S_1\big/K \cap S_1$ \ is semi-simple and is contained in the theme \ $E\big/K$. So it has rank \ $\leq 1$. As we know that \ $S_1$ \ is not contained in \ $K$, the rank is exactely \ $1$, and \ $K \cap S_1$ \ has corank \ $1$ \ in \ $S_1$. As \ $K \cap L_1 = \{0\}$ \ the previous proposition shows that \ $E\big/K \cap S_1$ \ is a rank \ $d$ \ theme. So the case 2  is proved.\\
For the proof of the case 3 it is enough to prove the equalities
\begin{align*}
&  S_1(E') = S_{j_0+1}(E)\big/L_{j_0} \quad L_1(E') = L_{j_0+1}\big/L_{j_0} 
\quad {\rm and} \\
& d(E') = d- j_0, \quad  rk(E') = rk(E) - j_0.
\end{align*}
As we know that \ $L_j$ \ is a theme of rank \ $j$, because \ $d(S_{j+1}) = j+1$, the rank of \ $E'$ \ is \ $ rk(E) - j_0$. The equality \ $d(L_{j_0}) = j_0$ \ is  then clear.\\
As \ $S_{j_0+1}(E)\big/L_{j_0} $ \ is semi-simple, we have \ $S_{j_0+1}(E)\big/L_{j_0} \subset S_1(E')$. But now we know that they have same rank because 
 $$rk(S_1(E')) = rk(E') - d(E')+1 = rk(E) - j_0 - (d-j_0) + 1 = rk(E) - d+ 1.$$
Also \ $L_1(E')$ \ and \ $L_{j_0+1}\big/L_{j_0}$ \ are rank 1  and normal submodules of \ $E'$. And as \ $S_2(E') = S_{j_0+2}\big/L_{j_0} $, we have \ $S_2(E')\big/\big[L_{j_0+1}\big/L_{j_0}\big] = S_{j_0+2}\big/S_{j_0+1}$ \ is semi-simple (rank 1), this gives the inclusion
$$L_1(E') : = \Sigma^1(S_2(E')) \subset L_{j_0+1}\big/L_{j_0} $$
and so the equality \ $L_1(E') =L_{j_0+1}\big/L_{j_0}$ \  is proved. $\hfill \blacksquare$

\bigskip

\begin{prop}[The rank 1 quotients]\label{rank 1 quot.}
Let \ $E$ \ be a \ $[\lambda]-$primitive rank \ $k$ \  fresco. Put \ $d : = d(E)$. Then any rank 1 quotient of \ $E$ \ is a rank 1 quotient of \ $E\big/\Sigma^1(E)$. As \ $E\big/\Sigma^1(E)$ \ is semi-simple of rank \ $k-d+1$ \ it shows that there are exactly \ $k-d+1$ \ isomorphism classes of such a rank 1 quotient and they corresponds to the fundamental invariants of \ $E\big/\Sigma^1(E)$ \ as follows : if \ $\lambda_1, \dots, \lambda_{k-d+1}$ \ are numbers associated to any  J-H. sequence of \ $E\big/\Sigma^1(E)$, the isomorphism classes of rank 1 quotients of \ $E$ \ are given by
$$ \lambda_1- k+d, \dots, \lambda_2 - k + d-1, \dots, \lambda_{k+d-1}.$$
\end{prop}

\parag{Proof} Let \ $H$ \ be a normal co-rank 1 submodule of \ $E$. As \ $E\big/H$ \ is semi-simple, \ $H$ \ contains \ $\Sigma^1(E)$ \ so \ $E\big/H$ \ is a rank 1 quotient of \ $E\big/\Sigma^1(E)$. The converse is obvious. $\hfill \blacksquare$

\parag{Remark} Assume \ $d(E) \geq 2$. With the exception of the (unique) rank 1 quotient of \ $E\big/S_1(E)$ \ which is \ $E\big/S_{d-1}(E)$, no rank 1 quotient of \ $E$ \ may be the rank 1 quotient of a quotient theme of rank \ $\geq 2$ \ of \ $E$. Another way to say that is the following : for any rank \ $r \geq 2$ \ quotient theme \ $T$ \ of \ $E$ \ we have \ $T\big/F_{r-1}(T) \simeq E\big/S_{d-1}(E)$.

\parag{Exemple} Let \ $E$ \ a rank \ $3$ \ $[\lambda]-$primitive fresco with \ $d(E) = 2$ \ (so \ $E$ \ is not semi-simple and is not a theme). Then there exists \ $k-d+1 = 2$ \ isomorphism classes of rank 1 quotients of \ $E$.\\
 For instance assume that  we have a J-H. sequence \ $0 \subset F_1 \subset F_2 \subset F_3 = E$, with \ $F_2 = S_1(E)$ \ and such \ $E\big/ F_1$ \ is a rank \ $2$ \ theme.\\
If \ $E \simeq (a -\lambda_1.b)S_1^{-1}.(a - \lambda_2.b)S_2^{-1}.(a - \lambda_3.b)$ \ this means, with \ $\lambda_{j+1} = \lambda_j + p_j -1$ \ for\ $ j = 1,2$, that \ $p_1 < 0$ \ or \ $p_1 \geq 1$ \ and no term in \ $b^{P_1}$ \ in \ $S_1$ \ and \ $p_2 \geq 0$ \ with a non zero term in \ $b^{p_2}$ \ in \ $S_2$.\\
So we  put \ $F_j\big/F_{j-1} \simeq E_{\lambda_j}$ \ for \ $j \in [1,3]$ ; then the rank 1 quotients are isomorphic to \ $E_{\lambda_3} \simeq E\big/S_1(E)$ \ or \ $E_{\lambda_1-2}$, because using the computations of section 3.3  we see that \ $\Sigma^1(E) \simeq E_{\lambda_2+1}$, and so \ $E_{\lambda_1-2}$ \ is a  rank 1 quotient of \ $E\big/\Sigma^1(E)$ \ and a fortiori of \ $E$.\\

To conclude we give a method to compute \ $L(E)$ \ is many cases.

\begin{lemma}\label{compute}
Let \ $E$ \ be a \ $[\lambda]-$primitive fresco and assume that \ $ E \simeq \A\big/\A.P$ \ where 
$$ P : = (a - \lambda_1.b).S_1^{-1} \dots S_{k-1}^{-1}.(a -\lambda_k.b) $$
where \ $\lambda_j + j$ \ is an increasing sequence. Assume that the first non commuting index is \ $h \in [1,k-1]$.  Then \ $E$ \ has a normal sub-theme of rank \ $2$ \ with fundamental invariants \ $\lambda_1, \lambda_{h+1}+ h$.
\end{lemma}

\parag{Proof} It is a simple application of the corollary \ref{commute 2} which also allow to compute the parameter of the rank \ $2$ \ obtain by commuting in \ $P$ \ from the parameter of the rank \ $2$ \ theme \ $F{h+1}\big/F_{h-1}$ \ and the integers \ $p_1, \dots, p_h$. $\hfill \blacksquare$

\parag{Remarks}
\begin{enumerate}[i)]
\item The hypothesis of the lemma means that \ $p_1, \dots ,p_{h-1}$ \ are positive and that for each \ $j \in [1,h-1]$ \ the coefficient of \ $b^{p_j}$ \ in \ $S_j$ \ is zero. But  the coefficient of \ $b^{p_h}$ \ in \ $S_h$ \ is not zero (and \ $p_h = 0$ \ is allowed).
\item The lemma implies that \ $L(E) \simeq E_{\lambda_1}$ \ is the first term of the principal J-H. sequence of \ $E$.
\end{enumerate}

   \section{Appendix : the existence theorem.}
   
   The aim of this appendix is to prove the the following existence theorem for  the fresco associated to a relative de Rham cohomology class :
   
   \begin{thm}\label{Existence}
 Let \ $X$ \ be a connected complex manifold of dimension \ $n + 1$ \ where \ $ n$ \ is a natural integer, and let \ $f : X \to D$ \ be an non constant proper  holomorphic function on an open  disc \ $D$ \ in \ $\mathbb{C}$ \ with center \ $0$. Let us assume that \ $df$ \ is nowhere vanishing outside of \ $X_0 : = f^{-1}(0)$.\\
 Let \ $\omega$ \ be a \ $\mathscr{C}^{\infty}-(p+1)-$differential form on \ $X$ \ such that \ $d\omega = 0 = df\wedge\omega $. Denote by \ $E$ \ the geometric \ $\A-$module \ $\mathbb{H}^{p+1}(X, (\hat{K}^{\bullet}, d^{\bullet}))$ \ and \ $[\omega]$ \ the image of \ $\omega$ \ in \ $E\big/B(E)$. 
 Then \ $\A.[\omega] \subset E\big/B(E)$ \ is a fresco.
 \end{thm}
 
 Note that this result is an obvious consequence of the finiteness theorem \ref{Finitude} that we shall prove below. It gives the fact that \ $E$ \ is naturally an \ $\A-$module which is of finite type over the subalgebra \ $\C[[b]]$ \ of \ $\A$, and so its \ $b-$torsion \ $B(E)$ \ is a finite dimensional \ $\C-$vector space. Moreover, the finiteness theorem asserts that \ $E \big/B(E)$ \ is a geometric (a,b)-module.

 \subsection{Preliminaries.}
 
 Here we shall complete and precise the results of the section 2 of [B.II]. The situation we shall consider is the following : let \ $X$ \ be a connected  complex manifold of dimension \ $n +1$ \ and \ $f : X \to \mathbb{C}$ \ a non constant holomorphic function such that \ $\{ x \in X / \ df = 0 \} \subset f^{-1}(0)$. We introduce the following complexes of sheaves supported by \ $X_0 : = f^{-1}(0)$
 \begin{enumerate}
 \item  The formal completion ''in \ $f$'' \ $(\hat{\Omega}^{\bullet}, d^{\bullet})$ \ of the usual holomorphic de Rham complex of \ $X$.
 \item The sub-complexes \ $(\hat{K}^{\bullet}, d^{\bullet})$ \ and \ $(\hat{I}^{\bullet}, d^{\bullet})$ \ of \  $(\hat{\Omega}^{\bullet}, d^{\bullet})$ \  where the subsheaves \ $\hat{K}^p$ \ and \ $\hat{I}^{p+1}$ \ are defined for each \ $p \in \mathbb{N}$ \  respectively as the kernel and the image of the map
 $$  \wedge df : \hat{\Omega}^p \to \hat{\Omega}^{p+1} $$
 given par exterior multiplication by \ $df$.  We have the exact sequence
 \begin{equation*} 0 \to (\hat{K}^{\bullet}, d^{\bullet}) \to (\hat{\Omega}^{\bullet},  d^{\bullet}) \to (\hat{I}^{\bullet}, d^{\bullet})[+1]  \to 0. \tag{1}
 \end{equation*}
 Note that \ $\hat{K}^0$ \ and \ $\hat{I}^0$ \ are zero by definition.
 \item The natural inclusions \ $\hat{I}^p \subset \hat{K}^p$ \ for all \ $p \geq 0$ \ are compatible with the diff\'erential \ $d$. This leads to an exact sequence of complexes
 \begin{equation*}
 0 \to (\hat{I}^{\bullet}, d^{\bullet}) \to (\hat{K}^{\bullet}, d^{\bullet}) \to ([\hat{K}/\hat{I}]^{\bullet}, d^{\bullet}) \to 0 .\tag{2}
 \end{equation*}
 \item We have a natural inclusion \ $f^*(\hat{\Omega}_{\mathbb{C}}^1) \subset \hat{K}^1\cap Ker\, d$, and this gives a sub-complex (with zero differential) of \ $(\hat{K}^{\bullet}, d^{\bullet})$. As in [B.07], we shall consider also the complex \ $(\tilde{K}^{\bullet}, d^{\bullet})$ \ quotient. So we have the exact sequence
 \begin{equation*}
  0 \to f^*(\hat{\Omega}_{\mathbb{C}}^1) \to (\hat{K}^{\bullet}, d^{\bullet}) \to (\tilde{K}^{\bullet}, d^{\bullet}) \to 0 . \tag{3}
  \end{equation*}
 We do not make the assumption here that \ $f = 0 $ \ is a reduced equation of \ $X_0$, and we do not assume that \ $n \geq 2$, so the cohomology sheaf in degree 1  of the complex \ $(\hat{K}^{\bullet}, d^{\bullet})$, which is equal to \ $\hat{K}^1 \cap Ker\, d$ \ does not co{\'i}ncide, in general with \ $f^*(\hat{\Omega}_{\mathbb{C}}^1)$. So the complex \ $ (\tilde{K}^{\bullet}, d^{\bullet})$ \ may have a non zero cohomology sheaf in degree 1.
  \end{enumerate} 
  Recall now that we have on the cohomology sheaves of the following complexes \\
   $(\hat{K}^{\bullet}, d^{\bullet}), (\hat{I}^{\bullet}, d^{\bullet}), ([\hat{K}/\hat{I}]^{\bullet}, d^{\bullet}) $ \ and \ $f^*(\hat{\Omega}_{\mathbb{C}}^1), (\tilde{K}^{\bullet}, d^{\bullet})$ \ natural operations \ $a$ \ and \ $b$ \ with the relation \ $a.b - b.a = b^2$. They are defined in a na{\"i}ve way by 
  $$  a : = \times f \quad {\rm and} \quad  b : = \wedge df \circ d^{-1} .$$
  The definition of \ $a$ \ makes sens obviously. Let me precise the definition of \ $b$ \ first in  the case of \ $\mathcal{H}^p(\hat{K}^{\bullet}, d^{\bullet})$ \ with \ $p \geq 2$  : if \ $x \in \hat{K}^p \cap Ker\, d$ \ write \ $x = d\xi$ \ with \ $\xi \in \hat{\Omega}^{p-1}$ \ and let \ $b[x] : = [df\wedge \xi]$. The reader will check easily that this makes sens.\\
  For \ $p = 1$ \ we shall choose \ $\xi \in \hat{\Omega}^0$ \ with the extra condition  that \ $\xi = 0$ \ on the smooth part of  \ $X_0$ \ (set theoretically). This is possible because the condition \ $df \wedge d\xi = 0 $ \ allows such a choice : near a smooth point of \ $X_0$ \ we can choose coordinnates such \ $ f = x_0^k$ \ and the condition on \ $\xi$ \ means independance of \ $x_1, \cdots, x_n$. Then \ $\xi$ \ has to be (set theoretically) locally constant on \ $X_0$ \ which is locally connected. So we may kill the value of such a \ $\xi$ \ along \ $X_0$.\\
  The case of the complex \ $(\hat{I}^{\bullet}, d^{\bullet})$ \ will be reduced to the previous one using the next  lemma.
  
  \begin{lemma}\label{tilde b}
  For each \ $p \geq 0$ \ there is a natural injective map
  $$\tilde{b} :  \mathcal{H}^p(\hat{K}^{\bullet}, d^{\bullet}) \to \mathcal{H}^p(\hat{I}^{\bullet}, d^{\bullet})$$
  which satisfies the relation \ $a.\tilde{b} = \tilde{b}.(b + a) $. For \ $p \not= 1$ \ this map is bijective.
  \end{lemma}
  
  \parag{Proof} Let \ $x \in \hat{K}^p \cap Ker\, d $ \ and write \ $x = d\xi $ \ where \ $x \in \hat{\Omega}^{p-1}$ \ (with \ $\xi = 0$ \ on \ $X_0$ \ if \ $p = 1$), and set \ $\tilde{b}([x]) : = [df\wedge \xi] \in \mathcal{H}^p(\hat{I}^{\bullet}, d^{\bullet})$. This is independant on the choice of \ $\xi$ \ because, for \ $p \geq 2$, adding \ $d\eta$ \ to \ $\xi$ \ does not modify the result as \ $[df\wedge d\eta] = 0 $. For \ $p =1$ \ remark that our choice of \ $\xi$ \ is unique.\\
   This is also independant of the the choice of \ $x $ \ in \ $ [x] \in \mathcal{H}^p(\hat{K}^{\bullet}, d^{\bullet})$ \  because adding \ $\theta \in \hat{K}^{p-1}$ \ to \ $\xi$ \ does not change \ $[df \wedge \xi]$.\\
   Assume \ $\tilde{b}([x]) = 0 $ \ in \ $ \mathcal{H}^p(\hat{I}^{\bullet}, d^{\bullet})$; this means that we may find \ $\alpha \in \hat{\Omega}^{p-2}$ \ such \ $df \wedge \xi = df \wedge d\alpha$. But then, \ $\xi - d\alpha $ \ lies in \ $\hat{K}^{p-1}$ \ and \ $x = d(\xi - d\alpha ) $ \ shows that \ $[x] = 0$. So \ $\tilde{b}$ \ is injective.\\
  Assume now \ $p \geq 2$.  If \ $df\wedge \eta $ \ is in \ $\hat{I}^p \cap Ker\, d$, then \ $df \wedge d\eta = 0 $ \ and \ $y : = d\eta $ \ lies in \ $\hat{K}^p \cap Ker\, d$ \ and defines a class \ $[y] \in  \mathcal{H}^p(\hat{K}^{\bullet}, d^{\bullet}) $ \ whose image by \ $\tilde{b}$ \ is \ $[df\wedge \eta] $. This shows the surjectivity of \ $\tilde{b}$ \ for \ $p \geq 2$.\\
   For \ $p=1$ \ the map \ $\tilde{b}$ \ is not surjective (see the remark below).\\
  To finish the proof let us  to compute \ $\tilde{b}(a[x] + b[x])$. Writing again \ $x = d\xi$, we get
   $$ a[x] + b[x] =[ f.d\xi + df \wedge \xi] = [d(f.\xi)] $$
   and so
   $$ \tilde{b}( a[x] + b[x] ) = [ df \wedge f.\xi ] = a.\tilde{b}([x]) $$
   which concludes the proof. $\hfill \blacksquare$
   
   \bigskip
   
   Denote by \ $i :  (\hat{I}^{\bullet}, d^{\bullet}) \to (\hat{K}^{\bullet}, d^{\bullet})$ \ the natural inclusion and define the action of \ $b$ \ on \ $\mathcal{H}^p(\hat{I}^{\bullet}, d^{\bullet})$ \ by \ $b : = \tilde{b}\circ \mathcal{H}^p(i) $. As \ $i$ \ is \ $a-$linear, we deduce the relation \ $a.b - b.a = b^2$ \ on \ $\mathcal{H}^p(\hat{I}^{\bullet}, d^{\bullet})$ \ from the relation of the previous lemma. \\
   
   The action of \ $a$ \ on the complex \ $ ([\hat{K}/\hat{I}]^{\bullet}, d^{\bullet}) $ \ is obvious and the action of \ $b$ \ is zero.\\
   
   The action of \ $a$ \ and \ $b$ \ on \ $f^*(\hat{\Omega}_{\mathbb{C}}^1) \simeq E_1\otimes \mathbb{C}_{X_0}$ \ are the obvious one, where \ $E_1$ \ is the rank 1 (a,b)-module with generator \ $e_1$ \ satisfying \ $a.e_1 = b.e_1$ \ (or, equivalentely, \ $E_1 : = \mathbb{C}[[z]]$ \ with \ $a : = \times z,\quad b : = \int_0^z $ \ and \ $e_1 : = 1$). \\
   Remark that the natural inclusion \ $f^*(\hat{\Omega}^1_{\mathbb{C}}) \hookrightarrow (\hat{K}^{\bullet}, d^{\bullet})$ \ is compatible with the actions of \ $a$ \ and \ $b$. The actions of \ $a$ \ and \ $b$ \ on \ $\mathcal{H}^1(\tilde{K}^{\bullet}, d^{\bullet}) $ \ are simply  induced by the corresponding actions on \ $\mathcal{H}^1(\hat{K}^{\bullet}, d^{\bullet})$.
   
   \parag{Remark} The exact sequence of complexes (1) induces  for any \ $p \geq 2$ \  a bijection
   $$ \partial^p : \mathcal{H}^p(\hat{I}^{\bullet}, d^{\bullet}) \to \mathcal{H}^p(\hat{K}^{\bullet}, d^{\bullet})$$
   and a short exact sequence 
   \begin{equation*}
    0 \to \mathbb{C}_{X_0} \to \mathcal{H}^1(\hat{I}^{\bullet}, d^{\bullet}) \overset{\partial^1}{\to} \mathcal{H}^1(\hat{K}^{\bullet}, d^{\bullet}) \to 0 \tag{@}
    \end{equation*}
   because of the de Rham lemma. Let us check  that for \ $p \geq 2$ \ we have \ $\partial^p = (\tilde{b})^{-1}$ \ and that for \ $p =1$ \ we have \ $\partial^1\circ \tilde{b} = Id$. If \ $x = d\xi \in \hat{K}^p \cap Ker\, d$ \ then \ $\tilde{b}([x]) = [df\wedge \xi]$ \ and \ $\partial^p[df\wedge\xi] = [d\xi]$. So \ $ \partial^p\circ\tilde{b} = Id \quad \forall p \geq 0$. For \ $p \geq 2$ \ and \ $df\wedge\alpha \in \hat{I}^p \cap Ker\, d$ \ we have \ $\partial^p[df\wedge\alpha] = [d\alpha]$ \ and \ $\tilde{b}[d\alpha] = [df\wedge\alpha]$, so \ $\tilde{b}\circ\partial^p = Id$. For \ $p = 1$ \ we have \ $\tilde{b}[d\alpha] = [df\wedge(\alpha - \alpha_0)]$ \ where \ $\alpha_0 \in \mathbb{C}$ \ is such that \ $\alpha_{\vert X_0} = \alpha_0$.  This shows that in degree 1 \ $\tilde{b}$ \ gives a canonical splitting of the exact sequence \ $(@)$.
   
  \subsection{$\A-$structures.}
   
   Let us consider now the \ $\mathbb{C}-$algebra
   $$ \A : = \{ \sum_{\nu \geq  0} \quad P_{\nu}(a).b^{\nu} \}$$
   where \ $P_{\nu} \in \mathbb{C}[z]$, and the commutation relation \ $a.b - b.a = b^2$, assuming that left and right multiplications by \ $a$ \ are continuous for the \ $b-$adic topology of \ $\A$.
   
   Define the following complexes of sheaves of left  \ $\A-$modules on \ $X$ :
   \begin{align*}
   & ({\Omega'}^{\bullet}[[b]], D^{\bullet})  \quad {\rm and} \quad   ({\Omega''}^{\bullet}[[b]], D^{\bullet}) \quad {\rm where} \tag{4}\\
    & {\Omega'}^{\bullet}[[b]] : = \sum_{j=0}^{+\infty} b^j.\omega_j \quad {\rm with} \quad \omega_0 \in \hat{K}^p  \\
     & {\Omega''}^{\bullet}[[b]] : = \sum_{j=0}^{+\infty} b^j.\omega_j \quad {\rm with} \quad \omega_0 \in \hat{I}^p \\
     & D(\sum_{j=0}^{+\infty} b^j.\omega_j) = \sum_{j=0}^{+\infty} b^j.(d\omega_j - df\wedge \omega_{j+1}) \\
     & a.\sum_{j=0}^{+\infty} b^j.\omega_j  = \sum_{j=0}^{+\infty} b^j.(f.\omega_j + (j-1).\omega_{j-1}) \quad {\rm with \ the \ convention} \quad \omega_{-1} = 0 \\
     & b.\sum_{j=0}^{+\infty} b^j.\omega_j  = \sum_{j=1}^{+\infty} b^j.\omega_{j -1}
     \end{align*}
     It is easy to check that \ $D$ \ is \ $\A-$linear and that \ $D^2 = 0 $. We have a natural inclusion of complexes of left \ $\A-$modules
     $$\tilde{i} :  ({\Omega''}^{\bullet}[[b]], D^{\bullet}) \to ({\Omega'}^{\bullet}[[b]], D^{\bullet}).$$
     
     Remark that we have natural morphisms of complexes
       \begin{align*}
     & u :   (\hat{I}^{\bullet}, d^{\bullet}) \to  ({\Omega''}^{\bullet}[[b]], D^{\bullet}) \\
     & v :  (\hat{K}^{\bullet}, d^{\bullet}) \to  ({\Omega'}^{\bullet}[[b]], D^{\bullet})
    \end{align*}
    and that these morphisms are compatible with \ $i$. More precisely, this means that we have the commutative diagram of complexes
 $$   \xymatrix{ (\hat{I}^{\bullet}, d^{\bullet}) \ar[d]^i \ar[r]^u & ({\Omega''}^{\bullet}[[b]], D^{\bullet}) \ar[d]^{\tilde{i}} \\
     (\hat{K}^{\bullet}, d^{\bullet})  \ar[r]^v &  ({\Omega'}^{\bullet}[[b]], D^{\bullet}) } $$
     
     The following theorem is a variant of theorem 2.2.1. of [B.II].

     \begin{thm}\label{(a,b)-structures}
     Let \ $X$ \ be a connected  complex manifold of dimension \ $n +1$ \ and \ $f : X \to \mathbb{C}$ \ a non constant holomorphic function such that 
      $$\{ x \in X / \ df = 0 \} \subset f^{-1}(0).$$
      Then the morphisms of complexes \ $u$ \ and \ $v$ \ introduced above are quasi-isomorphisms. Moreover, the isomorphims that they induce on the cohomology sheaves of these complexes are compatible with the actions of \ $a$ \ and \ $b$.
  \end{thm}
  
  \smallskip
  
  This theorem builds a natural structure of left \ $\A-$modules on each of the complex \\ 
  $(\hat{K}^{\bullet}, d^{\bullet}), (\hat{I}^{\bullet}, d^{\bullet}), ([\hat{K}/\hat{I}]^{\bullet}, d^{\bullet}) $ \ and \ $f^*(\hat{\Omega}_{\mathbb{C}}^1), (\tilde{K}^{\bullet}, d^{\bullet})$ \ in the derived category of bounded complexes of sheaves of \ $\mathbb{C}-$vector spaces on \ $X$.\\
   Moreover the short exact sequences
  \begin{align*}
&  0 \to (\hat{I}^{\bullet}, d^{\bullet}) \to (\hat{K}^{\bullet}, d^{\bullet}) \to ([\hat{K}/\hat{I}]^{\bullet}, d^{\bullet}) \to 0 \\
&  0 \to f^*(\hat{\Omega}_{\mathbb{C}}^1) \to (\hat{K}^{\bullet}, d^{\bullet}), (\hat{I}^{\bullet}, d^{\bullet}) \to (\tilde{K}^{\bullet}, d^{\bullet}) \to 0
\end{align*}
are equivalent to short exact sequences of complexes of left \ $\A-$modules in the derived category.

  \parag{Proof} We have to prove that for any \ $p \geq 0$ \ the maps \ $\mathcal{H}^p(u)$ \ and \ $\mathcal{H}^p(v)$ \ are bijective and compatible with the actions of \ $a$ \ and \ $b$. The case of \ $\mathcal{H}^p(v)$ \ is handled (at least for \ $n \geq 2$ \ and \ $f$ \ reduced) in prop. 2.3.1. of [B.II].  To seek for  completeness and for the convenience of the reader  we shall treat here the case of \ $\mathcal{H}^p(u)$.\\
  First we shall prove the injectivity of \ $\mathcal{H}^p(u)$. Let \ $\alpha = df\wedge \beta \in \hat{I}^p \cap Ker\, d$ \ and assume that we can find \ $U = \sum_{j=0}^{+\infty} b^j.u_j \in {\Omega''}^{p-1}[[b]]$ \ with \ $\alpha = DU$. Then we have the following relations
  $$ u_0 = df \wedge \zeta,  \quad \alpha = du_0 - df \wedge u_1 \quad {\rm and} \quad du_j = df\wedge u_{j+1} \ \forall j \geq 1. $$
  For \ $j \geq 1$ \ we have \ $[du_j] = b[du_{j+1}]$ \ in \ $ \mathcal{H}^p(\hat{K}^{\bullet}, d^{\bullet})$;  using corollary 2.2. of [B.II] which gives the \ $b-$separation of \ $\mathcal{H}^p(\hat{K}^{\bullet}, d^{\bullet})$, this implies \ $[du_j] = 0, \forall j \geq 1$ \ in \ $ \mathcal{H}^p(\hat{K}^{\bullet}, d^{\bullet})$. For instance we can find \ $\beta_1 \in \hat{K}^{p-1}$ \ such that \ $du_1 = d\beta_1$. Now, by de Rham, we can write \ $u_1 = \beta_1 + d\xi_1$ \ for \ $p \geq 2$, where \ $\xi_1 \in \hat{\Omega}^{p-2}$. Then we conclude that
   \ $ \alpha = -df\wedge d(\xi_1 + \zeta) $ \ and \ $[\alpha] = 0$ \ in \ $\mathcal{H}^p(\hat{I}^{\bullet}, d^{\bullet})$.\\
   For \ $p = 1$ \ we have \ $u_1 = 0 $ \ and \ $[\alpha] = [-df\wedge d\xi_1] = 0$ \ in \ $\mathcal{H}^1(\hat{I}^{\bullet}, d^{\bullet})$.\\
   We shall show now that the image of \ $\mathcal{H}^p(u)$ \ is dense in \ $\mathcal{H}^p({\Omega''}^{\bullet}[[b]], D^{\bullet})$ \  for its  \ $b-$adic topology. Let \ $\Omega : = \sum_{j=0}^{+\infty} \ b^j.\omega_j \in {\Omega''}^{p}[[b]]$ \ such that \ $D\Omega = 0$. The following relations holds \ $ d\omega_j = df\wedge \omega_{j+1} \quad \forall j \geq 0 $ \ and \ $\omega_0 \in \hat{I}^p$. The corollary 2.2. of [B.II] again allows to find \ $\beta_j \in \hat{K}^{p-1}$ \ for any \ $j \geq 0$ \ such that \ $d\omega_j = d\beta_j$. Fix \ $N \in \mathbb{N}^*$. We have
   $$ D(\sum_{j=0}^N b^j.\omega_j) = b^N.d\omega_N = D(b^N.\beta_N) $$
   and \ $\Omega_N : = \sum_{j=0}^N b^j.\omega_j  - b^N.\beta_N $ \ is \ $D-$closed and in \ ${\Omega''}^{p}[[b]]$. And we have \ $\Omega - \Omega_N \in b^N.\mathcal{H}^p({\Omega''}^{\bullet}[[b]], D^{\bullet})$, so the sequence \ $(\Omega_N)_{N \geq 1}$ \ converges to \ $\Omega$ \ in \ $\mathcal{H}^p({\Omega''}^{\bullet}[[b]], D^{\bullet})$ \ for its \ $b-$adic topology. Let us show that each \ $\Omega_N$ \ is in the image of \ $\mathcal{H}^p(u)$.\\
   Write \ $\Omega_N : = \sum_{j=0}^N b^j.w_j $. The condition \ $D\Omega_N = 0$ \ implies \ $dw_N = 0$ \ and \ $dw_{N-1} = df\wedge w_N = 0$. If we write \ $w_N = dv_N$ \ we obtain \ $d(w_{N-1} + df\wedge v_N) = 0$ \ and  \ $\Omega_N - D(b^N.v_N) $ \ is of degree \ $N-1$ \ in \ $b$. For \ $N =1$ \ we are left with \ $w_0 + b.w_1 - (-df\wedge v_1 + b.dv_1) = w_0 + df\wedge v_1$ \ which is in \ $\hat{I}^p \cap Ker\, d$ \ because \ $dw_0 = df\wedge dv_1$.\\
   To conclude it is enough to know the following two facts
   \begin{enumerate}[i)]
   \item The fact  that \ $\mathcal{H}^p(\hat{I}^{\bullet}, d^{\bullet})$ \ is complete for its \ $b-$adic topology.
   \item The fact that \ $Im(\mathcal{H}^p(u)) \cap b^N.\mathcal{H}^p({\Omega''}^{\bullet}[[b]], D^{\bullet}) \subset Im(\mathcal{H}^p(u)\circ b^N) \quad \forall N \geq 1 $.
   \end{enumerate}
   Let us first conclude the proof of the surjectivity of \ $\mathcal{H}^p(u)$ \ assuming i) and ii).\\
   For any \ $[\Omega] \in \mathcal{H}^p({\Omega''}^{\bullet}[[b]], D^{\bullet})$ \ we know that there exists a sequence \ $(\alpha_N)_{N \geq 1}$ \ in \ $ \mathcal{H}^p(\hat{I}^{\bullet}, d^{\bullet})$ \ with \ $\Omega - \mathcal{H}^p(u)(\alpha_N) \in b^N.\mathcal{H}^p({\Omega''}^{\bullet}[[b]], D^{\bullet})$. Now the property ii) implies that we may choose the sequence \ $(\alpha_N)_{N \geq 1} $ \ such that \ $[\alpha_{N+1}] - [\alpha_N]$ \ lies in \ $ b^N.\mathcal{H}^p(\hat{I}^{\bullet}, d^{\bullet})$. So the property i) implies that the Cauchy sequence \ $([\alpha_N])_{N \geq 1} $ \ converges to \ $[\alpha] \in \mathcal{H}^p(\hat{I}^{\bullet}, d^{\bullet})$. Then the continuity of \ $\mathcal{H}^p(u)$ \ for the \ $b-$adic topologies coming from its \ $b-$linearity, implies \ $\mathcal{H}^p(u)([\alpha]) = [\Omega]$.\\
   The compatibility with \ $a$ \ and \ $b$ \ of the maps \ $\mathcal{H}^p(u)$ \ and \ $\mathcal{H}^p(v)$ \ is an easy exercice. 
   
   \smallskip
   
   Let us now prove properties i) and ii).\\
   The property i) is a direct consequence of the completion of \  $\mathcal{H}^p(\hat{K}^{\bullet}, d^{\bullet})$ \ for its \ $b-$adic topology given by the corollary 2.2. of [B.II] \ and the \ $b-$linear isomorphism \ $\tilde{b} $ \ between \ $\mathcal{H}^p(\hat{K}^{\bullet}, d^{\bullet})$ \ and \ $\mathcal{H}^p(\hat{I}^{\bullet}, d^{\bullet})$ \ constructed in the lemma 2.1.1. above.\\
   To prove ii) let \ $\alpha \in \hat{I}^p \cap Ker\, d $ \ and \ $N \geq 1$ \ such that
   $$ \alpha = b^N.\Omega + DU $$
   where \ $\Omega \in {\Omega''}^{p}[[b]]$ \ satisfies \ $D\Omega = 0$ \ and where \ $U \in  {\Omega''}^{p-1}[[b]]$. With obvious notations we have
   \begin{align*}
   & \alpha = du_0 -df\wedge u_1\\
   & \cdots \\
   & 0  = du_j - df\wedge u_{j+1}  \quad \forall j \in [1, N-1] \\
   & \cdots \\
   & 0 = \omega_0 + du_N - df\wedge u_{N+1} 
   \end{align*}
   which implies \ $D(u_0+ b.u_1+ \cdots + b^N.u_N) = \alpha + b^N.du_N$ \ and the fact that \ $du_N$ \ lies in \ $\hat{I}^p \cap Ker \, d$. So we conclude that \ $[\alpha] + b^N.[du_N] $ \ is in the kernel of \ $\mathcal{H}^p(u)$ \ which is \ $0$. Then \ $[\alpha] \in b^N.\mathcal{H}^p(\hat{I}^{\bullet}, d^{\bullet})$.
   $\hfill \blacksquare$
   
     \parag{Remark} The map
  $$ \beta : ({\Omega'}[[b]]^{\bullet}, D^{\bullet}) \to ({\Omega''}[[b]]^{\bullet}, D^{\bullet})$$
  defined by \ $\beta(\Omega) = b.\Omega$ \ commutes to the differentials and with the action of \ $b$. It induces the isomorphism \ $\tilde{b}$ \ of the lemma \ref{tilde b} on the cohomology sheaves. So it is a quasi-isomorphism of complexes of \ $\mathbb{C}[[b]]-$modules.\\
 To prove this fact, it is enough to verify that the diagram\\
 $$ \xymatrix{&(\hat{K}^{\bullet}, d^{\bullet}) \ar[d]^{\tilde{b}} \ar[r]^v &  ({\Omega'}[[b]]^{\bullet}, D^{\bullet}) \ar[d]^{\beta} \\
 & (\hat{I}^{\bullet}, d^{\bullet}) \ar[r]^u &  ({\Omega''}[[b]]^{\bullet}, D^{\bullet})}$$
 
 \bigskip
 
 induces  commutative diagams on the cohomology sheaves. \\
  But this is clear because if \ $\alpha = d\xi$ \ lies in \ $ \hat{K}^p \cap Ker \, d$ \ we have \ $D(b.\xi) = b.d\xi - df\wedge \xi $ \ so \ $\mathcal{H}^p (\beta)\circ \mathcal{H}^p (v)([\alpha]) = \mathcal{H}^p (u)\circ \mathcal{H}^p (\tilde{b})([\alpha])$ \ in \ $\mathcal{H}^p ({\Omega''}[[b]]^{\bullet}, D^{\bullet}). \hfill \blacksquare$

    \subsection{The finiteness theorem.} 
  
  Let us recall some basic definitions on the left modules over the algebra \ $\A$.

  Now let \ $E$ \ be any left \ $\A-$module, and define \ $B(E)$ \ as the \ $b-$torsion of \ $E$. that is to say
  $$ B(E) : = \{ x \in E \ / \  \exists N \quad  b^N.x = 0 \}.$$
  Define \ $A(E)$ \ as the \ $a-$torsion of \ $E$ \ and 
   $$\hat{A}(E) : = \{x \in E \ / \ \mathbb{C}[[b]].x \subset A(E) \}.$$
   Remark that \ $B(E)$ \ and \ $\hat{A}(E)$ \ are  sub-$\A-$modules of \ $E$ \ but that \ $A(E)$ \ is not stable by \ $b$. 
  
  \begin{defn}\label{petit}
  A left \ $\A-$module \ $E$ \ is called {\bf small} when the following conditions hold
  \begin{enumerate}
  \item \ $E$ \ is a finite type \ $\mathbb{C}[[b]]-$module ;
  \item \ $B(E) \subset \hat{A}(E)$ ;
  \item \ $\exists N \ / \  a^N.\hat{A}(E) = 0 $ ;
  \end{enumerate}
  \end{defn}
  
  Recall that for \ $E$ \ small we have always the equality \ $ B(E) = \hat{A}(E)$ \ (see [B.I] lemme 2.1.2) and that this complex vector space is finite dimensional. The quotient \ $E/B(E)$ \ is an (a,b)-module called {\bf the associate (a,b)-module} to \ $E$.\\
  Conversely, any left \ $\A-$module \ $E$ \ such that \ $B(E)$ \ is a finite dimensional \ $\mathbb{C}-$vector space and such that \ $E/B(E)$ \ is an (a,b)-module is small.\\
  The following easy  criterium to be small will be used later :
  
  \begin{lemma}\label{crit. small}
  A left \ $\A-$module \ $E$ \ is small if and only if the following conditions hold :
  \begin{enumerate}
  \item \ $\exists N \ / \ a^N.\hat{A}(E) = 0 $ ;
  \item \ $B(E) \subset \hat{A}(E) $ ;
  \item  \ $\cap_{m\geq 0} b^m.E \subset \hat{A}(E) $ ;
  \item \ $Ker \, b$ \ and \ $Coker \, b$ \ are finite dimensional complex vector spaces.
  \end{enumerate}
  \end{lemma} 
  
  As the condition 3 in the previous lemma has been omitted in [B.II] (but this does not affect this article because this lemma was used only in a case were this condition 3 was satisfied, thanks to proposition 2.2.1. of {\it loc. cit.}), we shall give the (easy) proof.
  \parag{Proof} First the conditions 1 to 4  are obviously necessary. Conversely, assume that \ $E$ \ satisfies these four conditions. Then condition 2 implies that the action of \ $b$ \ on \ $\hat{A}(E)\big/B(E)$ \ is injective. But the condition 1 implies that \ $b^{2N} = 0$ \ on \ $\hat{A}(E) $ \ (see [B.I] ). So we conclude that \ $\hat{A}(E) = B(E) \subset Ker\, b^{2N}$ \ which is a finite dimensional complex vector space using condition 4 and an easy induction. Now \ $E/B(E)$ \ is a \ $\mathbb{C}[[b]]-$module which is separated for its \ $b-$adic topology. The finitness of \ $Coker \, b$ \ now shows that it is a free finite type \ $\mathbb{C}[[b]]-$module concluding the proof. $\hfill \blacksquare$
  
  \begin{defn}\label{geometric}
  We shall say that a left \ $\A-$module \ $E$ \ is {\bf geometric} when \ $E$ \ is small and when it associated (a,b)-module \ $E/B(E)$ \ is geometric.
  \end{defn}
  
  The main result of this section is the following theorem, which shows that the Gauss-Manin connection of a proper holomorphic function produces geometric \ $\A-$modules associated to vanishing cycles and nearby cycles.

 \begin{thm}\label{Finitude}
 Let \ $X$ \ be a connected complex manifold of dimension \ $n + 1$ \ where \ $ n$ \ is a natural integer, and let \ $f : X \to D$ \ be an non constant proper  holomorphic function on an open  disc \ $D$ \ in \ $\mathbb{C}$ \ with center \ $0$. Let us assume that \ $df$ \ is nowhere vanishing outside of \ $X_0 : = f^{-1}(0)$.\\
 Then the \ $\A-$modules 
 $$ \mathbb{H}^j(X, (\hat{K}^{\bullet}, d^{\bullet})) \quad {\rm and} \quad \mathbb{H}^j(X, (\hat{I}^{\bullet}, d^{\bullet})) $$
 are geometric for any \ $j \geq 0 $.
 \end{thm}
 
In the proof we shall use the \ $\mathscr{C}^{\infty}$ \ version of the complex \ $(\hat{K}^{\bullet}, d^{\bullet})$. We define \ $K_{\infty}^p$ \ as the kernel of \ $\wedge df : \mathscr{C}^{\infty,p} \to \mathscr{C}^{\infty,p+1}$ \ where \ $\mathscr{C}^{\infty,j}$ \ denote the sheaf of \ $\mathscr{C}^{\infty}-$ \ forms on \ $X$ \ of degree p, let \ $\hat{K}_{\infty}^p$ \ be the \ $f-$completion and \ $(\hat{K}_{\infty}^{\bullet}, d^{\bullet})$ \ the corresponding de Rham complex.

  The next lemma is proved in [B.II] (lemma  6.1.1.)
 
 \begin{lemma}\label{diff}
 The natural inclusion
 $$ (\hat{K}^{\bullet}, d^{\bullet}) \hookrightarrow (\hat{K}_{\infty}^{\bullet}, d^{\bullet}) $$
 induce a quasi-isomorphism.
 \end{lemma}
 
 \parag{Remark} As the sheaves \ $\hat{K}_{\infty}^{\bullet}$ \ are fine, we have a natural isomorphism
 $$ \mathbb{H}^p(X, (\hat{K}^{\bullet}, d^{\bullet})) \simeq H^p\big(\Gamma(X, \hat{K}_{\infty}^{\bullet}), d^{\bullet}\big).$$
 
 Let us denote by \ $X_1$ \ the generic fiber of \ $f$. Then \ $X_1$ \ is a smooth compact complex manifold of dimension \ $n$ \ and the restriction of \ $f$ \ to \ $f^{-1}(D^*)$ \ is a locally trivial \ $\mathscr{C}^{\infty}$ \ bundle with typical fiber \ $X_1$ \ on \ $D^* = D \setminus \{0\}$, if the disc \ $D$ \ is small enough around \ $0$. Fix now \ $\gamma \in H_p(X_1, \mathbb{C})$ \ and let \ $(\gamma_s)_{s \in D^*}$ \ the corresponding multivalued horizontal family of \ $p-$cycles \ $\gamma_s \in H_p(X_s, \mathbb{C})$. Then, for \ $\omega \in \Gamma(X, \hat{K}_{\infty}^p \cap Ker\, d)$, define the multivalued holomorphic function
 $$ F_{\omega}(s) : = \int_{\gamma_s} \frac{\omega}{df} .$$
 Let now 
  $$\Xi : = \oplus_{\alpha \in \mathbb{Q} \cap ]-1,0], j \in [0,n]} \quad  \mathbb{C}[[s]].s^{\alpha}.\frac{(Log s)^j}{j!} .$$
  This is an \ $\A-$modules with \ $a$ \ acting as multiplication by \ $s$ \ and \ $b$ \ as the primitive in \ $s$ \ without constant. Now if \ $\hat{F}_{\omega}$ \ is the asymptotic expansion at \ $0$ \ of \ $F_{\omega}$, it is an element  in \ $\Xi$, and we obtain in this way an \ $\A-$linear map
 $$ Int :   \mathbb{H}^p(X, (\hat{K}^{\bullet}, d^{\bullet})) \to H^p(X_1, \mathbb{C}) \otimes_{\mathbb{C}} \Xi .$$
 To simplify notations, let  \ $E : =  \mathbb{H}^p(X, (\hat{K}^{\bullet}, d^{\bullet}))$. Now using Grothendieck theorem [G.65], there exists \ $N \in \mathbb{N}$ \ such that \ $Int(\omega) \equiv 0 $, implies \ $a^N.[\omega] = 0$ \ in \ $E$.  As the converse is clear we conclude that \ $\hat{A}(E) =  Ker(Int)$. It is also clear that \ $B(E) \subset Ker(Int)$ \ because \ $\Xi$ \ has no \ $b-$torsion. So we conclude that \ $E$ \ satisfies properties 1 and 2 of the lemma \ref{crit. small}. The property 3 is also true because of the regularity of the Gauss-Manin connection of \ $f$.
 
 \parag{End of the proof of theorem \ref{Finitude}} To show that \ $E : =  \mathbb{H}^p(X, (\hat{K}^{\bullet}, d^{\bullet}))$ \ is small, it is enough to prove that \ $E$ \ satisfies the condition 4 of the lemma \ref{crit. small}. Consider now the long exact sequence of hypercohomology of the exact sequence of complexes
 $$   0 \to (\hat{I}^{\bullet}, d^{\bullet}) \to (\hat{K}^{\bullet}, d^{\bullet}) \to ([\hat{K}/\hat{I}]^{\bullet}, d^{\bullet}) \to 0 .$$
 It contains the exact sequence
 $$  \mathbb{H}^{p-1}(X, ([\hat{K}\big/\hat{I}]^{\bullet}, d^{\bullet})) \to \mathbb{H}^p(X, (\hat{I}^{\bullet}, d^{\bullet})) \overset{\mathbb{H}^p(i)}{\to} \mathbb{H}^p(X, (\hat{K}^{\bullet}, d^{\bullet})) \to \mathbb{H}^{p}(X, ([\hat{K}\big/\hat{I}]^{\bullet}, d^{\bullet})) $$
 and we know that \ $b$ \ is induced on the complex of \ $\A-$modules  quasi-isomorphic to \ $(\hat{K}^{\bullet}, d^{\bullet})$ \ by the composition \ $i\circ \tilde{b}$ \ where \ $\tilde{b}$ \ is a quasi-isomorphism of complexes of \ $\mathbb{C}[[b]]-$modules. This implies that the kernel  and the cokernel of \ $\mathbb{H}^p(i)$ \ are isomorphic (as \ $\mathbb{C}-$vector spaces) to \ $Ker\, b$ \ and \ $Coker \, b$ \ respectively. Now to prove that \ $E$ \ satisfies condition 4 of the lemma \ref{crit. small} it is enough to prove finite dimensionality for the vector spaces \ $  \mathbb{H}^{j}(X, ([\hat{K}\big/\hat{I}]^{\bullet}, d^{\bullet})) $ \ for all \ $j \geq 0 $.\\
 But the sheaves \ $[\hat{K}\big/\hat{I}]^j \simeq [Ker\,df\big/Im\, df]^j$ \ are coherent on \ $X$ \ and supported in \ $X_0$. The spectral sequence
 $$ E_2^{p,q} : = H^q\big( H^p(X, [\hat{K}\big/\hat{I}]^{\bullet}), d^{\bullet}\big) $$ 
 which converges to \ $ \mathbb{H}^{j}(X, ([\hat{K}\big/\hat{I}]^{\bullet}, d^{\bullet})) $, is a bounded complex of finite dimensional vector spaces by Cartan-Serre. This gives the desired finite dimensionality.\\
 To conclude the proof, we want to show that \ $E/B(E)$ \ is geometric. But this is an easy consequence of the regularity of the Gauss-Manin connexion of \ $f$ \ and of the Monodromy theorem, which are already incoded in the definition of \ $\Xi$ : the injectivity on \ $E/B(E)$ \ of the \ $\A-linear$ \ map \ $Int$ \ implies that \ $E/B(E)$ \ is geometric. \\
 Remark now that the piece of  exact sequence above gives also the fact that \ $\mathbb{H}^p(X, (\hat{I}^{\bullet}, d^{\bullet}))$ \ is geometric, because it is an exact sequence of \ $\A-$modules. $\hfill \blacksquare$
 
 \parag{Remark} It is easy to see that the properness assumption on \ $f$ \ is only used for two purposes : \\
 --To have a (global) \ $\mathscr{C}^{\infty}$ \ Milnor fibration on a small punctured disc around \ $0$, with a finite dimensional cohomology for the Milnor fiber.\\
 -- To have compactness of the singular set \ $\{df = 0 \}$.\\
 This allows to give with the same proof an analoguous finiteness result in many other situations.
  
\bigskip

\section*{Bibliography}

\begin{itemize}

\item{[Br.70]} Brieskorn, E. {\it Die Monodromie der Isolierten Singularit{\"a}ten von Hyperfl{\"a}chen}, Manuscripta Math. 2 (1970), p. 103-161.\\

\item{[B.84]}  Barlet, D. \textit{Contribution du cup-produit de la fibre de Milnor aux p\^oles de \ $\vert f \vert^{2\lambda}$},  Ann. Inst. Fourier (Grenoble) t. 34, fasc. 4 (1984), p. 75-107.\\

\item{[B. 93]} Barlet, D. {\it Th\'eorie des (a,b)-modules I}, in Complex Analysis and Geo-metry, Plenum Press, (1993), p. 1-43.\\

\item{[B. 95]} Barlet, D. {\it Th\'eorie des (a,b)-modules II. Extensions}, in Complex Analysis and Geometry, Pitman Research Notes in Mathematics Series 366 Longman (1997), p. 19-59.\\

\item{[B. 05]} Barlet, D. {\it Module de Brieskorn et forme hermitiennes pour une singularit\'e isol\'ee d'hypersuface}, revue de l'Inst. E. Cartan (Nancy) 18 (2005), p. 19-46. \\

\item{[B.I]} Barlet, D. {\it Sur certaines singularit\'es non isol\'ees d'hypersurfaces I}, Bull. Soc. math. France 134 (2), 2006, p.173-200.\\

\item{[B. II]} Barlet, D. {\it Sur certaines singularit\'es d'hypersurfaces II}, J. Alg. Geom. 17 (2008), p. 199-254.

\item{[B. 08]} Barlet, D. {\it Two finiteness theorem for regular (a,b)-modules}, preprint Institut E. Cartan (Nancy) (2008) $n^0 5$, p. 1-38, arXiv:0801.4320 (math. AG and math. CV)

\item{[B.09]} Barlet,D. {\it P\'eriodes \'evanescentes et (a,b)-modules monog\`enes}, Bollettino U.M.I. (9) II (2009) p.651-697.\\

\item{[B.III]} Barlet, D. {\it Sur les fonctions \`a lieu singulier de dimension 1 }, Bull. Soc. math. France 137 (4), 2009, p.587-612.\\

\item{[B.10]} Barlet,D. {\it Le th\`eme d'une p\'eriode \'evanescente}, preprint Institut E. Cartan (Nancy) (2009) $n^0 33$, p.1-57.\\

\item{[B.-S. 04]} Barlet, D. et Saito, M. {\it Brieskorn modules and Gauss-Manin systems for non isolated hypersurface singularities,} J. Lond. Math. Soc. (2) 76 (2007) $n^01$ \ p. 211-224.\\

\item{[Bj.93]} Bj{\"o}rk, J-E, {\it   Analytic D-modules and applications}, Kluwer Academic  publishers (1993).\\

 \item{[K.76]} Kashiwara, M. {\it b-function and holonomic systems}, Inv. Math. 38 (1976) p. 33-53.\\

\item{[M.74]} Malgrange, B. {\it Int\'egrale asymptotique et monodromie}, Ann. Sc. Ec. Norm. Sup. 7 (1974), p.405-430.\\

\item{[M. 75]} Malgrange, B. {\it Le polyn\^ome de Bernstein d'une singularit\'e isol\'ee}, in Lect. Notes in Math. 459, Springer (1975), p.98-119.\\

\item{[S. 89]} Saito, M. {\it On the structure of Brieskorn lattices}, Ann. Inst. Fourier 39 (1989), p.27-72.\\

\end{itemize}

\end{document}